\documentclass[a4paper]{amsart}

\usepackage[ngerman,french,american]{babel}
\usepackage[T1]{fontenc}
\usepackage[latin1]{inputenc}
\usepackage{amssymb}
\usepackage[all]{xy}
\usepackage{amsxtra}
\usepackage{enumerate}

\usepackage[lite]{amsrefs}

\newcommand*{\brd}{-\hspace{0pt}}
\newcommand*{\nbd}{\nobreakdash-\hspace{0pt}}

\newcommand*{\N}{{\mathbb{N}}}

\newcommand*{\R}{{\mathbb{R}}}
\newcommand*{\C}{{\mathbb{C}}}

\newcommand*{\Cat}{\mathcal{C}}

\newcommand*{\braket}[2]{\langle#1,#2\rangle}

\newcommand*{\inv}{^\times}
\newcommand*{\blank}{{\llcorner\!\!\lrcorner}}

\newcommand*{\ID}{{\mathrm{id}}}
\newcommand*{\Gl}{{\mathrm{Gl}}}
\newcommand*{\SC}{{\mathrm{SC}}}
\newcommand*{\op}{{\mathrm{op}}}
\newcommand*{\modular}{{\mu}}

\newcommand*{\CCINF}{\mathcal{D}}
\newcommand*{\CINF}{\mathcal{E}}
\newcommand*{\CONT}{\mathcal{C}}
\newcommand*{\CCONT}{\mathcal{C}_{\mathrm{c}}}
\newcommand*{\UniEnvel}{\mathcal{U}}
\newcommand*{\LG}{\mathfrak{g}}
\newcommand*{\Mult}{\mathcal{M}}
\newcommand*{\LMult}{\mathcal{M}_{\mathrm{l}}}
\newcommand*{\RMult}{\mathcal{M}_{\mathrm{r}}}

\newcommand*{\RepG}{\mathbf{R}_G}
\newcommand*{\NRepG}{\hat{\mathbf{R}}_G}
\newcommand*{\RepH}{\mathbf{R}_H}
\newcommand*{\NRepH}{\hat{\mathbf{R}}_H}
\newcommand*{\EssG}{\mathbf{M}_G}
\newcommand*{\EssH}{\mathbf{M}_H}
\newcommand*{\ModG}{\hat{\mathbf{M}}_G}
\newcommand*{\ModH}{\hat{\mathbf{M}}_H}

\newcommand*{\IN}{\smallint\!\!}

\newcommand*{\hot}{\mathbin{\hat{\otimes}}}
\newcommand*{\prot}{\mathbin{\hat{\otimes}_\pi}}
\newcommand*{\rprod}{\prod\nolimits'}
\newcommand*{\defeq}{\mathrel{:=}}

\newcommand*{\prto}{\twoheadrightarrow}
\newcommand*{\into}{\rightarrowtail}

\DeclareMathOperator{\supp}{supp}

\DeclareMathOperator{\Map}{Map}
\DeclareMathOperator{\Ker}{Ker}
\DeclareMathOperator{\Coker}{Coker}
\DeclareMathOperator{\Hom}{Hom}
\DeclareMathOperator{\End}{End}
\DeclareMathOperator{\Aut}{Aut}
\DeclareMathOperator{\smooth}{S}
\DeclareMathOperator{\rough}{R}
\DeclareMathOperator{\PREC}{Pt}
\DeclareMathOperator{\Center}{Z}
\DeclareMathOperator{\ci}{Ic}
\DeclareMathOperator{\nci}{I}
\DeclareMathOperator{\Ind}{Ind}
\DeclareMathOperator{\cInd}{c-Ind}
\DeclareMathOperator{\Res}{Res}
\DeclareMathOperator{\Tor}{Tor}
\DeclareMathOperator{\Ext}{Ext}
\DeclareMathOperator{\Ho}{H}
\DeclareMathOperator{\Left}{\mathbb{L}}
\DeclareMathOperator{\Right}{\mathbb{R}}

\theoremstyle{plain}
\newtheorem{theorem}{Theorem}[section]
\newtheorem{proposition}[theorem]{Proposition}
\newtheorem{lemma}[theorem]{Lemma}
\newtheorem{corollary}[theorem]{Corollary}

\theoremstyle{definition}
\newtheorem{definition}[theorem]{Definition}
\newtheorem{remark}[theorem]{Remark}

\begin{document}

\title[Smooth group representations]{Smooth group representations on
  bornological vector spaces}

\author{Ralf Meyer}
 
\address{\selectlanguage{ngerman}
         Mathematisches Institut\\
         Westfälische Wilhelms-Universität Münster\\
         Einsteinstraße 62\\
         48149 Münster\\
         Germany}

\email{rameyer@math.uni-muenster.de}

\thanks{This research was supported by the EU-Network \emph{Quantum
  Spaces and Noncommutative Geometry} (Contract HPRN-CT-2002-00280)
  and the \emph{\selectlanguage{ngerman}Deutsche
  Forschungsgemeinschaft} (SFB 478).}

\begin{abstract}
  We develop the basic theory of smooth representations of locally
  compact groups on bornological vector spaces.  In this setup, we are
  able to formulate better general theorems than in the topological
  case.  Nonetheless, smooth representations of totally disconnected
  groups on vector spaces and of Lie groups on Fréchet spaces remain
  special cases of our theory.  We identify smooth representations
  with essential modules over an appropriate convolution algebra.  We
  examine smoothening functors on representations and modules and show
  that they agree if they are both defined.  We establish the basic
  properties of induction and compact induction functors using adjoint
  functor techniques.  We describe the center of the category of
  smooth representations.

  \medskip

  {\selectlanguage{french} Nous développons la théorie basique des
  représentations lisses des groupes localement compacts sur les
  espaces vectorielles bornologiques.  Dans ce contexte, nous pouvons
  établir des meilleurs théorèmes que dans la situation topologique.
  Néanmoins, les représentations lisses des groupes totalement
  discontinus sur les espaces vectorielles et les représentations
  lisses des groupes de Lie sur les espaces de Fréchet restent des cas
  spéciales de notre théorie.  Nous identifions des représentations
  lisses avec des modules essentielles sur une algèbre de convolution
  convenable.  Nous examinons des foncteurs régularisants sur des
  représentations et des modules et nous montrons que ils sont égales
  si ils sont définis.  Nous établons les propriétés basiques des
  foncteurs d'induction et d'induction compact en employant des
  techniques des foncteurs adjointes.  Nous décrivons le centre du
  catégorie des représentations lisses.}
\end{abstract}

\subjclass[2000]{22D12, 22D20, 46A17, 18H10}

\keywords{center, multiplier algebra, induced representation,
  smoothening, bornological metrizability}

\maketitle

\section{Introduction}
\label{sec:introduction}

Smooth representations of totally disconnected groups on vector spaces
and of Lie groups on locally convex topological vector spaces have
already been studied for a long time.  It is also known that one can
define smooth representations of arbitrary locally compact groups
using the spaces of smooth functions introduced by François Bruhat
in~\cite{Bruhat:Distributions}.  We shall consider, instead, smooth
representations of locally compact groups on
\emph{bornological} vector spaces
(see~\cite{Hogbe-Nlend:Bornologies}).  While this may appear to be
only a minor variation on the usual theory, it turns out that there
are several small but significant details that make the bornological
theory much more pleasant and more powerful.  Smooth representations
of totally disconnected groups on vector spaces and of Lie groups on
Fréchet spaces are special cases of our theory, so that it allows for
a unified treatment of these two kinds of representations.

Bornological vector spaces went out of fashion quite some time ago.
This is rather unfortunate because they are the ideal setting for
noncommutative geometry.  As soon as we move beyond Fréchet spaces, we
run into annoying problems when we work with topological vector
spaces.  For instance, the multiplication on an algebra like
$\CCINF(\R)$ with convolution is only separately continuous and not
jointly continuous.  Therefore, one has to give \emph{ad hoc}
definitions for the complexes that compute the Hochschild and cyclic
homology of such convolution algebras.  Problems of this nature are
artefacts which disappear if we work bornologically instead.
Moreover, bornologies are essential for the purposes of local cyclic
cohomology, which is a variant of cyclic cohomology that produces
better results for Banach algebras like the algebra of continuous
functions on a compact space.

A great advantage of bornological versus topological analysis is the
adjoint associativity between the completed bornological tensor
product~$\hot$ and the internal $\Hom$ functor: $\Hom(A\hot
B,C)\cong\Hom\bigl(A,\Hom(B,C)\bigr)$.  In particular, there is a
canonical bornology on the space $\Hom(B,C)$ of bounded linear maps
between two bornological vector spaces.  Adjoint associativity holds
for vector spaces and Banach spaces, but not for topological vector
spaces.  It provides bornological analysis with a much richer
algebraic structure than topological analysis.  For representation
theory this means that the general theory of smooth representations
of locally compact groups on bornological vector spaces is very
similar to the purely algebraic theory of smooth representations of
totally disconnected groups on vector spaces.

An instance of this is our main theorem, which asserts that the
category of smooth representations of~$G$ is isomorphic to the
category of essential modules over the convolution algebra $\CCINF(G)$
of smooth functions with compact support on~$G$.  We also have very
nice adjointness relations between restriction, induction and compact
induction functors, from which we can deduce many properties of these
functors.

We now explain our results in greater detail.  Throughout this
article, $G$ denotes a locally compact topological group.  Bruhat
(\cite{Bruhat:Distributions}) defines spaces $\CCINF(G)$ and
$\CINF(G)$ of smooth functions with compact support and with arbitrary
growth at infinity, respectively.  In the totally disconnected case a
function is smooth if and only if it is locally constant.  In the Lie
group case smoothness has the usual meaning.  General locally compact
groups are treated using the deep structure theory of almost connected
groups.  We recall Bruhat's definitions and adapt them to our
bornological setup in Section~\ref{sec:smooth_functions_on_groups}.
Besides basic facts about these function spaces, we prove some
interesting results about metrizable bornological vector spaces.

A representation $\pi\colon G\to\Aut(V)$ on a complete convex
bornological vector space~$V$ is called smooth if the map that sends
$v\in V$ to the function $g\mapsto \pi(g,v)$ takes values in
$\CINF(G,V)$ and is a bounded linear map $\pi_\ast\colon
V\to\CINF(G,V)$.  Equivalently, the formula $Wf(g)\defeq g\cdot f(g)$
defines a bounded linear operator on $\CCINF(G,V)$.  For totally
disconnected~$G$ this amounts to the requirement that any bounded set
be stabilized by an open subgroup of~$G$.  In particular, if~$V$ is a
vector space with the fine bornology, we get the usual notion of a
smooth representation of a totally disconnected group on a complex
vector space.

Now suppose~$G$ to be a Lie group.  A representation is called
differentiable if it is $k$~times continuously differentiable for all
$k\in\N$.  This notion is weaker than smoothness.  For instance, the
left regular representation on the space of compactly supported
distributions $\CINF'(G)$ is differentiable but not smooth.
Differentiability and smoothness are equivalent if~$V$ is
bornologically metrizable.  In particular, this happens if~$V$ is a
Fréchet space equipped with a reasonable bornology.

Differentiable representations on bornological vector spaces are
closely related to smooth representations on topological vector
spaces.  We show that a bornological representation~$\pi$ is
differentiable if and only if it extends to a bounded algebra
homomorphism $\IN\pi\colon \CINF'(G)\to\End(V)$.  Similarly, a
topological representation~$\pi$ is smooth if and only if it extends
to a bounded homomorphism $\IN\pi\colon \CINF'(G)\to\End(V)$, where
$\End(V)$ carries the equicontinuous bornology.  Let~$V$ be a
bornological topological vector space, equip it with the von Neumann
bornology.  Then there is no difference between the spaces of
continuous and bounded maps $V\to V$, equipped with the equicontinuous
and equibounded bornology, respectively.  Hence topological smoothness
is equivalent to bornological differentiability in this case.
If~$V$ is a Fréchet space, we know that bornological differentiability
and smoothness are equivalent, so that the topological and
bornological notions of smooth representation agree for Fréchet
spaces.  For general~$V$ the bornological notion of smoothness is
more restrictive than the topological one.

If we restrict~$\IN\pi$ to the convolution algebra $\CCINF(G)$, we
turn~$V$ into a module over $\CCINF(G)$.  A module~$V$ over
$\CCINF(G)$ is called \emph{essential} if the module action is a
bornological quotient map $\CCINF(G)\hot V\to V$.  That is, each
bounded subset of~$V$ is the image of a bounded subset of
$\CCINF(G)\hot V$.  The following theorem generalizes a well-known and
much used fact for totally disconnected groups:

\begin{theorem}  \label{the:essential_smooth_intro}
  Let~$G$ be a locally compact group.  The categories of essential
  bornological left $\CCINF(G)$\brd{}modules and of smooth
  representations of~$G$ on bornological vector spaces are isomorphic.
  The isomorphism sends a representation $\pi\colon G\to\Aut(V)$ to
  the module $\IN\pi\colon\CCINF(G)\to\End(V)$.
\end{theorem}

The theorem makes three assertions.  First, if $\pi\colon G\to\Aut(V)$
is smooth, then $\IN\pi\colon \CCINF(G)\hot V\to V$ is a bornological
quotient map.  In fact, this map even has a bounded linear section.
Secondly, any essential module over $\CCINF(G)$ arises in this fashion
from a smooth representation of~$G$.  Thirdly, a bounded linear map
between two smooth representations is $G$\nbd{}equivariant if and only
if it is a homomorphism of $\CCINF(G)$\brd{}modules.  In the topological
framework it is still true that~$\pi$ is smooth if and only if
$\IN\pi\colon \CCINF(G,V)\to V$ has a continuous linear section
(see~\cite{Blanc:Homologie}) .  However, $\CCINF(G,V)$ is no longer a
topological tensor product of $\CCINF(G)$ and~$V$.  Therefore, we fail
to characterize smooth representations in terms of the algebra
$\CCINF(G)$.

We study analogues in the category of modules over $\CCINF(G)$ of
several constructions with representations, namely, smoothening,
restriction, induction and compact induction.  Let $H\subseteq G$ be
a closed subgroup.  Then we only have $\CCINF(H)\subseteq\CINF'(G)$,
so that the restriction of a $\CCINF(G)$\brd{}module to a
$\CCINF(H)$\brd{}module is not always defined.  If~$V$ is an arbitrary
$\CCINF(G)$\brd{}module, then $\CCINF(G)\hot_{\CCINF(G)} V$ and
$\Hom_{\CCINF(G)}(\CCINF(G),V)$ carry canonical $\CCINF(H)$\brd{}module
structures.  The resulting functors are called the \emph{smooth and
  rough restriction} functors, $\smooth_G^H$ and $\rough_G^H$.  In the
converse direction, if~$V$ is a module over $\CCINF(H)$, we can
produce a module over $\CCINF(G)$ in two ways.  We define the
\emph{compact induction functor} and the \emph{rough induction
  functor} by
\begin{align*}
  \ci_H^G(V) &\defeq \CCINF(G)\hot_{\CCINF(H)} V.
  \\
  \nci_H^G(V) &\defeq \Hom_{\CCINF(H)}(\CCINF(G),V).
\end{align*}
The functors $\smooth\defeq\ci_G^G=\smooth_G^G$ and
$\rough\defeq\nci_G^G=\rough_G^G$ are called \emph{smoothening} and
\emph{roughening}, respectively.  Up to a relative modular factor,
$\smooth\circ \nci_H^G$ and $\ci_H^G$ agree with the induction and
compact induction functors for representations, respectively.

These functors enjoy many useful algebraic properties.  For instance,
they are exact for appropriate classes of extensions.  The exactness
of the smoothening functor implies that the class of essential modules
is closed under extensions.  The content of the roughening functor is
the following: roughly speaking, the roughening of a module~$V$ is the
largest module~$W$ that satisfies $\smooth V=\smooth W$.  Many
important properties of the induction and restriction functors follow
easily by playing around with adjoint associativity.  We prove the
Shapiro Lemma in group homology and cohomology in this fashion and we
show how to reduce $\Tor$ and $\Ext$ for the category of essential
$\CCINF(G)$\brd{}modules to group homology and cohomology.  It is
remarkable that such results can be proved easily and purely
algebraically.  There are no analytical difficulties whatsoever.

The smoothening functors for representations and modules also agree.
The module smoothening is the range of the map $\IN\pi\colon
\CCINF(G)\hot V\to V$.  The image of the uncompleted tensor product is
known as the Gårding subspace of~$V$.  Jacques Dixmier and Paul
Malliavin show in~\cite{Dixmier-Malliavin:Factorisations} that the
Gårding subspace is equal to the smoothening for Lie group
representations on Fréchet spaces.  The same is true for arbitrary
continuous representations of locally compact groups on bornological
vector spaces.

Finally, we examine the analogue of the Bernstein center of a totally
disconnected group.  This is the center of the category of smooth
representations of~$G$ on complex vector spaces, which was studied
first by Joseph Bernstein (\cite{Bernstein:Center}).  It plays a
crucial role in the representation theory of reductive groups, which
is parallel to the role played by the center of the universal
enveloping algebra in the Lie group case.

We prove that the center of the category of smooth representations
of~$G$ is isomorphic to the center of the multiplier algebra of
$\CCINF(G)$.  In the totally disconnected case this is the same as
the Bernstein center.  We describe the multiplier algebra of
$\CCINF(G)$ and its center as spaces of distributions on~$G$.  For Lie
groups the multiplier algebra is just $\CINF'(G)$.  For a connected
complex Lie group with trivial center, central multipliers are
necessarily supported at the identity element.  Thus the center of the
category of smooth representations of~$G$ is isomorphic to the center
of the universal enveloping algebra of~$G$ in this case.

\section{Spaces of smooth functions on locally compact groups}
\label{sec:smooth_functions_on_groups}

Many results of this section are adaptations to the bornological
setting of results of François Bruhat (\cite{Bruhat:Distributions}).
There are a few issues regarding tensor products and metrizability that
do not arise in the topological setting, however.

Since we are only dealing with complete convex bornologies, we drop
these adjectives from our notation: whenever we assert or ask that a
space be a bornological vector space, it is understood that it is
asserted or asked to be a complete convex bornological vector space.
Good references for the basic theory of bornological vector spaces are
the publications of Henri Hogbe-Nlend (\cites{Hogbe-Nlend:Completions,
Hogbe-Nlend:Theorie, Hogbe-Nlend:Bornologies}), whose notation we will
follow mostly.

\subsection{Preliminaries}
\label{sec:preliminaries}

The structure theory of locally compact groups is crucial for Bruhat's
definitions in order to reduce to the case of Lie groups.  Although
its results are very difficult to prove, they are extremely simple to
apply and state.

Let~$G$ be a locally compact group.  Let $G_0\subseteq G$ be the
connected component of the identity element.  The group~$G$ is called
\emph{totally disconnected} if $G_0=\{1\}$, \emph{connected} if
$G_0=G$ and \emph{almost connected} if $G/G_0$ is compact.

A totally disconnected locally compact group has a base for the
neighborhoods of the identity element consisting of compact open
subgroups (see~\cite{Hewitt-Ross:Harmonic_Analysis}).  Applying this
to the totally disconnected group $G/G_0$, we find that any locally
compact group contains an almost connected open subgroup.

\begin{theorem}[\cite{Montgomery-Zippin}]
  \label{the:structure_almost_connected}
  Let~$G$ be an almost connected locally compact group.  Then~$G$ is
  isomorphic to a projective limit of Lie groups.  More explicitly,
  there is a directed set~$I$ of compact normal subgroups $k\subseteq
  G$ such that $G/k$ is a Lie group for all $k\in I$ and $\bigcap I =
  \{1\}$.  We have $G=\varprojlim_{k\in I} G/k$ for any such system.
\end{theorem}

\begin{definition}  \label{def:smooth_subgroup}
  A subgroup $k\subseteq G$ is called \emph{smooth} if its normalizer
  $N_G(k)\subseteq G$ is open and $N_G(k)/k$ is a Lie group.  Let
  $\SC$ or $\SC(G)$ be the set of all smooth compact subgroups.  A
  \emph{fundamental system of smooth compact subgroups} in~$G$ is a
  set~$I$ of smooth compact subgroups which is directed by inclusion
  and satisfies $\bigcap I=\{1\}$.
\end{definition}

\begin{lemma}  \label{lem:smooth_subgroup}
  Let~$G$ be a locally compact group.  If $k\subseteq G$ is a smooth
  subgroup, then $G/k$, $k\backslash G$ and $G//k\defeq k\backslash
  G/k$ are smooth manifolds in a canonical way.  If $k_1\subseteq
  k_2$, then the induced maps $G/k_1\to G/k_2$, etc., are smooth.

  The set $\SC(G)$ is a fundamental system of smooth compact
  subgroups and in particular directed.  We have
  $$
  G\cong
  \varprojlim G/k \cong
  \varprojlim k\backslash G \cong
  \varprojlim G//k,
  $$
  where the limits are taken for $k\in\SC(G)$.

  A set of subgroups is a fundamental system of smooth compact
  subgroups if and only if it is a cofinal subset of $\SC(G)$.  The
  set~$I$ can be taken countable and even a decreasing sequence if and
  only if~$G$ is metrizable.
\end{lemma}

\begin{proof}
  Let $k\subseteq G$ be a smooth subgroup and let~$U$ be its
  normalizer.  Thus~$U$ is an open subgroup of~$G$, $k$ is a normal
  subgroup of~$U$ and $U/k$ is a Lie group.  The homogeneous space
  $G/k$ is just a disjoint union of copies $gU/k$ of the Lie group
  $U/k$ for $g\in G/U$ and hence a smooth manifold.  The same applies
  to $k\backslash G$.  The proof of the corresponding assertion for
  $G//k$ is more complicated.  We view this as the orbit space of the
  action of~$k$ on $G/k$ by left multiplication.  For any $g\in G$,
  let $k'\defeq k\cap gUg^{-1}$.  Then $k\backslash kgU/k\cong
  k'\backslash gU/k\cong g^{-1}k'g\backslash U/k$ because $G/U$ is
  open.  The latter double coset space is really a left coset space
  because~$k$ is normal in~$U$.  Thus $k\backslash G/k$ is a disjoint
  union of smooth manifolds as well.
  
  Let $U\subseteq G$ be an open almost connected subgroup.  For~$U$
  instead of~$G$, our assertions follow from
  Theorem~\ref{the:structure_almost_connected}.  Since
  $\SC(U)\subseteq\SC(G)$ is cofinal, the latter is a fundamental
  system of smooth compact subgroups in~$G$.  We also get the
  isomorphisms $G\cong \varprojlim G/k$, etc., from the corresponding
  statement for~$U$.  It is clear that any cofinal subset of $\SC(G)$
  is still a fundamental system of smooth compact subgroups.
  Conversely, if~$I$ is such a set, then $I\subseteq\SC(G)$.  Let
  $k\in\SC(G)$.  Since $\cap I=\{1\}$, the set of $k'\in I$ with
  $k'\subseteq N_G(k)$ is cofinal. Since the Lie group $N_G(k)/k$ does
  not contain arbitrarily small subgroups, the quotient group $k'/k$
  must eventually be trivial, that is, $k'\subseteq k$.  This means
  that~$I$ is cofinal in $\SC(G)$.  It is clear from
  $G\cong\varprojlim G/k$ that~$G$ is metrizable if and only if we can
  choose~$I$ countable.
\end{proof}

Before we can define smooth functions on locally compact groups, we
need some generalities about spaces of smooth functions on manifolds
(see~\cite{Meyer:Born_Top} for more details).  Let~$M$ be a smooth
manifold and let~$B$ be a Banach space.  Then we equip the space
$\CCINF(M,B)$ of smooth functions with compact support from~$M$ to~$B$
with the following bornology.  A set~$S$ of smooth functions is
bounded if all $f\in S$ are supported in a fixed compact subset of~$M$
and the set of functions $D(S)$ is uniformly bounded for any
differential operator~$D$ on~$M$.  This is the von Neumann bornology
for the usual LF\brd{}topology on $\CCINF(M,B)$.  We let $\CCINF(M)$
be $\CCINF(M,\R)$ or $\CCINF(M,\C)$, depending on whether we work with
real or complex bornological vector spaces.  In the following, we will
assume that we work with complex vector spaces, but everything works
for real vector spaces as well.

If~$V$ is a bornological vector space, we let $\CCINF(M,V)$ be the
space of all functions $M\to V$ that belong to $\CCINF(M,V_T)$ for
some bounded complete disk $T\subseteq V$.  A subset of $\CCINF(M,V)$
is bounded if it is bounded in $\CCINF(M,V_T)$ for some~$T$.  (Recall
that~$V_T$ is the linear span of~$T$ equipped with the norm whose
closed unit ball is~$T$.  Hence it is a Banach space.)

Let~$\hot$ be the completed projective bornological tensor product.
It is defined by the universal property that bounded linear maps
$A\hot B\to C$ correspond to bounded bilinear maps $A\times B\to C$.
The natural map $\CCINF(M)\hot B\to\CCINF(M,B)$ is a bornological
isomorphism for all Banach spaces~$B$.  The functor
$\CCINF(M)\hot\blank$ commutes with direct limits and preserves
injectivity of linear maps because $\CCINF(M)$ is nuclear
(see~\cite{Hogbe-Nlend-Moscatelli:Nuclear}).  Therefore, we have
\begin{equation}  \label{eq:CCINF_M_tensor}
  \CCINF(M,V) \cong \CCINF(M)\hot V
\end{equation}
for all bornological vector spaces~$V$.  Moreover, for two
manifolds $M_1,M_2$ we have
$$
\CCINF(M_1)\hot\CCINF(M_2)\cong \CCINF(M_1\times M_2).
$$

We define the spaces $\CCONT^k(M,V)$ of $k$~times continuously
differentiable functions with compact support similarly for $k\in\N$.
If~$V$ is a Banach space, we let $\CCONT^k(M,V)$ be the usual LF-space
and equip it with the von Neumann bornology.  For general~$V$ we let
$\CCONT^k(M,V)\defeq \varinjlim \CCONT^k(M,V_T)$.  We let
$\CCONT^\infty(M,V)\defeq \varprojlim \CCONT^k(M,V)$ and call
functions in $\CCONT^\infty(M,V)$ \emph{differentiable} (see
also~\cite{Waelbroeck:Differentiable}).  While there evidently is no
difference between smooth functions and $\CONT^\infty$\brd{}functions
with values in a Banach space, smoothness is more restrictive than
differentiability in general.  Smooth functions are easier to work
with because of~\eqref{eq:CCINF_M_tensor}, which fails for
$\CCONT^\infty(M,V)$.

\begin{definition}  \label{def:metrizable}
  A bornological vector space is \emph{metrizable} if for any sequence
  $(S_n)$ of bounded subsets there is a sequence of scalars
  $(\epsilon_n)$ such that $\sum \epsilon_n S_n$ is bounded.
\end{definition}

The precompact bornology and the von Neumann bornology on a Fréchet
space are metrizable in this sense (see~\cite{Meyer:Born_Top}).

\begin{lemma}  \label{lem:differentiable_metrizable}
  If~$V$ is metrizable, then $\CCINF(M,V)=\CCONT^\infty(M,V)$.
\end{lemma}

\begin{proof}
  Let $S\subseteq\CCONT^\infty(M,V)$ be bounded.  That is, $S$ is
  bounded in $\CCONT^k(M,V)$ for all $k\in\N$.  For any $k\in\N$,
  there is a bounded complete disk $T_k\subseteq V$ such that~$S$ is
  bounded in $\CCONT^k(M,V_{T_k})$.  By metrizability, we can absorb
  all~$T_k$ in some bounded complete disk $T\subseteq V$.  Thus~$S$
  is bounded in $\CCONT^k(M,V_T)$ for all $k\in\N$.  This means
  that~$S$ is bounded in $\CCINF(M,V)$.
\end{proof}

\begin{lemma}  \label{lem:metrizable_product}
  A bornological vector space~$V$ is metrizable if and only if the
  functor $V\hot\blank$ commutes with countable direct products.
\end{lemma}

\begin{proof}
  It is easy to see that~$V$ is metrizable once $V\hot \prod_\N
  \C\cong \prod_\N (V\hot\C)$.  For the converse implication, we
  clearly have a bounded linear map $V\hot \prod B_n\to \prod V\hot
  B_n$.  We have to show that $\prod V\hot B_n$ satisfies the
  universal property of $V\hot \prod B_n$.  That is, we need that a
  bounded bilinear map $l\colon V\times \prod B_n\to X$ induces a
  bounded linear map $\prod V\hot B_n\to X$.  By definition, a bounded
  subset~$S$ of $\prod_\N V\hot B_n$ is contained in $\prod S_n\hot
  T_n$ with bounded complete disks $S_n$ and~$T_n$ in $V$ and~$B_n$,
  respectively.  Here $S_n\hot T_n$ denotes the complete disked hull
  of $S_n\times T_n$ in $V\hot B_n$.  By metrizability, all~$S_n$ are
  absorbed by some bounded complete disk $S'\subseteq V$.  Moving
  the absorbing constants into~$T_n$, we obtain $S\subseteq S'\hot
  \prod T_n'$.  This implies the desired universal property.
\end{proof}

\subsection{The definitions of the function spaces}
\label{sec:CCINF}

Let~$G$ be a locally compact group and let~$V$ be a bornological
vector space.  The spaces $\CCINF(G/k,V)$ are defined for all
$k\in\SC(G)$.  We pull back functions on $G/k$ to~$G$ and thus view
$\CCINF(G/k,V)$ as a space of functions on~$G$.  If $k_1\subseteq
k_2$, then $\CCINF(G/k_2,V)$ is the subspace of
right-$k_2$\brd{}invariant functions in $\CCINF(G/k_1,V)$ and thus a
retract of $\CCINF(G/k_1,V)$.  The set $\SC$ is directed by
Lemma~\ref{lem:smooth_subgroup}.  Hence the spaces $\CCINF(G/k,V)$ for
$k\in\SC$ form a strict inductive system.  Strict means that the
structure maps are bornological embeddings.  We let $\CCINF(G,V)$ be
its inductive limit.  This is just the union of the spaces
$\CCINF(G/k,V)$ equipped with the direct union bornology and thus a
space of $V$\nbd{}valued functions on~$G$.  We get the same space if we
replace $\SC$ by any fundamental system of smooth compact subgroups
because the latter are cofinal subsets of $\SC$.  In particular,
if~$G$ is metrizable, then we can use a decreasing sequence of
subgroups.

\begin{lemma}  \label{lem:CCINF_symmetric}
  We have
  $$
  \CCINF(G,V)=
  \varinjlim \CCINF(G/k,V) =
  \varinjlim \CCINF(k\backslash G,V) =
  \varinjlim \CCINF(G//k,V).
  $$
\end{lemma}

\begin{proof}
  For any compact subset $S\subseteq G/k$ there is $k_2\in\SC$ that
  stabilizes all points of~$S$.  That is, functions in $\CCINF(G/k,V)$
  with support in~$S$ are automatically left-$k_2$\brd{}invariant and
  hence belong to $\CCINF(G//k_2,V)$.  This yields the assertions.
\end{proof}

Let $H\subseteq G$ be a closed subgroup.  We define $\CCINF(G/H,V)$
and $\CCINF(H\backslash G,V)$ as follows.  The double coset space
$k\backslash G/H$ can be decomposed as a disjoint union of homogeneous
spaces for Lie groups as in the proof of
Lemma~\ref{lem:smooth_subgroup} and hence is a smooth manifold for all
$k\in I$.  We view $\CCINF(k\backslash G/H)$ as a space of
left-$k$\brd{}invariant functions on $G/H$.  If $k_1\subseteq k_2$,
then $\CCINF(k_2\backslash G/H)$ is the set of
left-$k_2$\brd{}invariant functions in $\CCINF(k_1\backslash G/H)$.
Thus the spaces $\CCINF(k\backslash G/H)$ for $k\in I$ form a strict
inductive system.  We let $\CCINF(G/H,V)\defeq \varinjlim
\CCINF(k\backslash G/H,V)$.  The definition of $\CCINF(H\backslash
G,V)$ is analogous.  Lemma~\ref{lem:CCINF_symmetric} shows that this
reproduces the old definition of $\CCINF(G/H,V)$ if~$H$ is normal
in~$G$.  If~$H$ is a compact subgroup, then $\CCINF(G/H,V)$ is
canonically isomorphic to the space $\CCINF(G,V)^H$ of elements in
$\CCINF(G,V)$ that are invariant under right translation by~$H$.

If~$G$ is a Lie group, then $G/H$ is a smooth manifold and
$\CCINF(G/H,V)$ evidently agrees with the usual space of smooth
functions defined in Section~\ref{sec:preliminaries}.  If~$G$ is
totally disconnected, then the spaces $k\backslash G/H$ are discrete.
Therefore, $\CCINF(G/H,V)$ is the space of locally constant functions
with compact support from $G/H$ to~$V$.

\begin{definition}  \label{def:CINF}
  A function $f\colon G/H\to V$ is called smooth if $h\cdot
  f\in\CCINF(G/H,V)$ for all $h\in\CCINF(G/H)$.  We let $\CINF(G/H,V)$
  be the space of smooth functions from $G/H$ to~$V$.  A subset~$S$ of
  $\CINF(G/H,V)$ is bounded if $h\cdot S$ is bounded in $\CCINF(G/H,V)$
  for all $h\in\CCINF(G/H)$.  We let $\CINF(G/H)\defeq\CINF(G/H,\C)$.
  
  For a closed subset $S\subseteq G/H$, let $\CINF_0(S,V)$ be the
  subspace of $\CINF(G/H,V)$ of functions supported in~$S$ and let
  $\CINF(S,V)$ be the quotient of $\CINF(G/H,V)$ by the ideal of
  functions vanishing in~$S$.  (The latter notation is slightly
  ambiguous because $\CINF(S,V)$ also depends on $G/H$.)
\end{definition}

Let $S\subseteq G/H$ be compact.  Then there is $h\in\CCINF(G/H)$ with
$h|_S=1$.  Therefore, we obtain the same spaces $\CINF_0(S,V)$ and
$\CINF(S,V)$ if we replace $\CINF(G,V)$ by $\CCINF(G,V)$ in the above
definition.  It is evident that $\CCINF(G/H,V)=\varinjlim
\CINF_0(S,V)$ where~$S$ runs through the directed set of compact
subsets of $G/H$.  Thus $\CCINF(G/H,V)$ is the space of compactly
supported elements of $\CINF(G/H,V)$.  However, the space
$\CINF(G/H,V)$ tends to be harder to analyze than $\CCINF(G/H,V)$.

\subsection{Nuclearity and exactness properties}
\label{sec:CCINF_nuclear}

Next we examine some properties of $\CCINF(G/H)$ and of the functor
$V\mapsto\CCINF(G/H,V)$.  Since the bornological tensor product
commutes with inductive limits, \eqref{eq:CCINF_M_tensor} implies
\begin{equation}
  \label{eq:CCINF_tensor}
  \CCINF(G/H,V) \cong
  \varinjlim \CCINF(k\backslash G/H,V) \cong
  \varinjlim \CCINF(k\backslash G/H)\hot V \cong
  \CCINF(G/H) \hot V.
\end{equation}

\begin{proposition}  \label{pro:CCINF_nuclear}
  The bornological vector space $\CCINF(G/H)$ is nuclear.  More
  generally, if~$V$ is nuclear, so is $\CCINF(G/H,V)$.
\end{proposition}

\begin{proof}
  For $k\in\SC$ and $S\subseteq k\backslash G/H$ compact, the subspace
  $\CINF_0(S)\subseteq \CCINF(k\backslash G/H)$ is a nuclear Fréchet
  space because $k\backslash G/H$ is a smooth manifold.  Hence it is
  nuclear as a bornological vector space as well
  (see~\cite{Hogbe-Nlend-Moscatelli:Nuclear}).  As an inductive limit
  of these spaces, the space $\CCINF(G/H)$ is nuclear as well.  Since
  nuclearity is hereditary for tensor products,
  \eqref{eq:CCINF_tensor} implies that $\CCINF(G/H,V)$ is nuclear
  if~$V$ is.
\end{proof}

To state the exactness properties of the functor $\CCINF(G/H,\blank)$,
we recall some natural classes of extensions.  A \emph{bornological
extension} is a diagram $K\overset{i}\to E\overset{p}\to Q$ with
$i=\Ker p$ and $p=\Coker i$.  It is called \emph{linearly split} if it
has a bounded linear section.  Then it follows that $E\cong K\oplus
Q$.  It is called \emph{locally linearly split} if for any bounded
complete disk $T\subseteq Q$ there is a local bounded linear section
$Q_T\to E$ defined on the Banach space~$Q_T$.  Equivalently, the
sequence
$$
0\to\Hom(B,K)\to\Hom(B,E)\to\Hom(B,Q)\to0
$$
is exact for any Banach space~$B$.  Locally linearly split
extensions are important for local cyclic cohomology.

\begin{proposition}  \label{pro:CCINF_exact}
  The functor $V\mapsto \CCINF(G/H,V)$ commutes with direct limits.
  It preserves bornological extensions and injectivity of morphisms.
  It also preserves locally linearly split and linearly split
  extensions.
\end{proposition}

\begin{proof}
  For any bornological vector space~$W$, the functor $V\mapsto W\hot
  V$ commutes with direct limits and preserves linearly split and
  locally linearly split bornological extensions.  Nuclearity of~$W$
  implies that it also preserves injectivity of morphisms and
  bornological extensions.  This yields the assertions because
  of~\eqref{eq:CCINF_tensor}.
\end{proof}

Now we turn from $\CCINF(G/H,V)$ to $\CINF(G/H,V)$.  For any open
covering of~$G/H$ there is a subordinate partition of unity consisting
of functions in $\CCINF(G/H)$.  In order to avoid taking square roots,
our convention for partitions of unity is that $\sum \phi_j^2(x)=1$.
We choose such a partition of unity $(\phi_j)_{j\in J}$ on $G/H$ with
$\phi_j\in\CCINF(G/H)$ for all $j\in J$ and use it to define maps
\begin{equation}
  \label{eq:CINF_retract}
  \begin{alignedat}{2}
    \iota &\colon \CINF(G/H,V)\to\prod_{j\in J} \CCINF(G/H,V),
    \qquad & \iota(f)_j &\defeq f\cdot\phi_j,
    \\
    \pi &\colon \prod_{j\in J} \CCINF(G/H,V)\to\CINF(G/H,V),
    \qquad & \pi\bigl((f_j)\bigr) &\defeq \sum_{j\in J} f_j\cdot\phi_j.
  \end{alignedat}
\end{equation}
It is clear that~$\iota$ is a well-defined bounded linear map.  The
map~$\pi$ is a well-defined bounded linear map as well because all but
finitely many of the products $f_j\phi_j h$ vanish for
$h\in\CCINF(G/H)$.  Thus $\CINF(G/H,V)$ is naturally isomorphic to a
retract (that is, direct summand) of $\prod_{j\in J} \CCINF(G/H,V)$.

\begin{proposition}  \label{pro:CINF_exact}
  The functor $\CINF(G/H,\blank)$ preserves bornological extensions
  and injectivity of morphisms.  It also preserves locally linearly
  split and linearly split bornological extensions.  The space
  $\CINF(G/H,V)$ is nuclear if (and only if) $V$ is nuclear and $G/H$
  is countable at infinity.
\end{proposition}

\begin{proof}
  The classes of extensions that occur in the proposition are closed
  under direct products.  Hence a retract of a direct product of exact
  functors is again exact.  Using the maps in~\eqref{eq:CINF_retract},
  the assertions about $\CINF(G/H,\blank)$ therefore follow from the
  corresponding assertions about $\CCINF(G/H,\blank)$ in
  Proposition~\ref{pro:CCINF_exact}.  Suppose $G/H$ to be countable at
  infinity.  Then the partition of unity above is countable, so that
  $\CINF(G/H,V)$ is a retract of a countable direct product of spaces
  $\CCINF(G/H,V)$.  Since nuclearity is hereditary for
  \emph{countable} direct products, $\CINF(G/H,V)$ is nuclear.
\end{proof}

\begin{definition}  \label{def:compact_support}
  Let $l\colon \CCINF(G/H,V)\to W$ be a bounded linear map.  Its
  \emph{support} $\supp l$ is the smallest closed subset $S\subseteq
  G/H$ such that $l(f)=0$ for all $f\in\CCINF(G/H,V)$ that vanish in a
  neighborhood of~$S$.  (An argument using partitions of unity shows
  that this is well-defined.)
  
  Let $\CCINF'(G/H,V)$ be the dual space of $\CCINF(G/H,V)$, equipped
  with the equibounded bornology.  For $S\subseteq G/H$, let
  $\CCINF'(S,V)\subseteq\CCINF'(G/H,V)$ be the set of linear
  functionals supported in~$S$.  Let $\CINF'(G/H,V)\defeq \varinjlim
  \CCINF'(S,V)$, where~$S$ runs through the compact subsets of $G/H$.
  In particular, for $V=\C$, we obtain the spaces $\CCINF'(G/H)$ and
  $\CINF'(G/H)$ of \emph{distributions} and \emph{distributions with
  compact support} on~$G/H$.
\end{definition}

\begin{lemma}  \label{lem:CINF_dual}
  The natural map from the dual of $\CINF(G/H,V)$ to $\CCINF'(G/H,V)$
  is a bornological isomorphism onto $\CINF'(G/H,V)$.  In particular,
  $\CINF'(G/H)$ is the dual space of $\CINF(G/H)$.
\end{lemma}

\begin{proof}
  It is not hard to see that for any set of bornological vector spaces
  $(V_x)$, the dual space of $\prod V_x$ is bornologically isomorphic
  to the direct sum $\bigoplus V_x'$.  This together
  with~\eqref{eq:CINF_retract} yields the assertion.
\end{proof}

\begin{lemma}  \label{lem:CINF_metrizable}
  If $G/H$ is countable at infinity and~$V$ is metrizable, then
  $$
  \CINF(G/H,V)\cong \CINF(G/H)\hot V.
  $$
\end{lemma}

\begin{proof}
  We have already shown that $\CCINF(G/H,V)\cong\CCINF(G/H)\hot V$.
  Using the maps in~\eqref{eq:CINF_retract} and
  Lemma~\ref{lem:metrizable_product}, we obtain $\CINF(G/H,V)\cong
  \CINF(G/H)\hot V$ as well.
\end{proof}

However, $\CINF(G/H,V)$ is not isomorphic to $\CINF(G/H)\hot V$ in
general.  All three spaces $\CINF(G/H\times G/H)$,
$\CINF(G/H,\CINF(G/H))$ and $\CINF(G/H)\hot\CINF(G/H)$ are different
unless $G/H$ is a smooth manifold or compact.  This is the reason why
the regular representation on $G/H$ usually fails to be smooth.

\subsection{Functoriality with respect to the group}
\label{sec:CCINF_functoriality}

\begin{definition}  \label{def:smooth_map}
  A continuous linear map $f\colon G_1/H_1\to G_2/H_2$ between two
  homogeneous spaces is called \emph{smooth} if for any $x\in G_1/H_1$
  and any $k_2\in\SC(G_2)$, there is $k_1\in\SC(G_1)$ and an open
  $k_1$\brd{}invariant neighborhood $V\subseteq G_1/H_1$ of~$x$ such
  that the restriction of~$f$ to~$V$ descends to a smooth map
  $k_1\backslash V\to k_2\backslash G_2/H_2$.
\end{definition}

\begin{lemma}  \label{lem:CCINF_functoriality}
  A smooth map $f\colon G_1/H_1\to G_2/H_2$ induces a bounded linear
  map
  $$
  f^\ast\colon \CINF(G_2/H_2,V)\to\CINF(G_1/H_1,V),
  \qquad f^\ast(h) \defeq h\circ f.
  $$
  If~$f$ is proper as well, $f^\ast$ restricts to a bounded linear map
  $$
  f^\ast\colon \CCINF(G_2/H_2,V)\to\CCINF(G_1/H_1,V),
  \qquad f^\ast(h) \defeq h\circ f.
  $$
\end{lemma}

\begin{proof}
  Use smooth partitions of unity.
\end{proof}

The following are examples of smooth maps.  They induce maps on spaces
of smooth functions by Lemma~\ref{lem:CCINF_functoriality}.

\begin{enumerate}
\item The group multiplication is a smooth map $G\times G\to G$.  So
  are the multiplication maps $G\times G/H\cong G\times G/1\times H\to
  G/H$ and $H\backslash G\times G\cong H\times 1\backslash G\times
  G\to H\backslash G$.  The map $G\times G\to G\times G$,
  $(x,y)\mapsto (x,xy)$, is smooth and so are similar maps involving
  homogeneous spaces.

\item The inversion is a smooth map $G\to G$ and $G/H\leftrightarrow
  H\backslash G$.

\item Any continuous group homomorphism is smooth.

\item If $g^{-1}Hg\subseteq H'$, then the map $G/H\to G/H'$ that sends
  $xH$ to $xHgH'=xgH'$ is smooth.

\end{enumerate}

Thus we can define the left and right \emph{regular representations}
$\lambda$ and~$\rho$ of~$G$ on $\CCINF(G,V)$ and $\CINF(G,V)$ by
\begin{equation}  \label{eq:biregular_representation}
  \lambda_g f(x)\defeq f(g^{-1}\cdot x),
  \qquad
  \rho_g f(x) \defeq f(x \cdot g).
\end{equation}

\begin{lemma}  \label{lem:CINF_invariant}
  The space $\CINF(G/H,V)$ is naturally isomorphic to the subspace of
  $\CINF(G,V)$ of functions~$f$ that satisfy $\rho_hf=f$ for all $h\in
  H$.
\end{lemma}

\begin{proof}
  The projection $G/H\to G$ is smooth and therefore induces a bounded
  injective map $\CINF(G/H,V)\to\CINF(G,V)$, whose range clearly
  consists of right-$H$\brd{}invariant functions.  Let $k\in\SC$ and
  let~$U$ be its normalizer.  In order to prove that $k\backslash G/H$
  is a smooth manifold, we decomposed $k\backslash G/H$ into a
  disjoint union of the double coset spaces $k\backslash UgH/H$ for
  $g\in U\backslash G/H$ and identified the contribution of each
  double coset with a homogeneous space for a Lie group action.  This
  reduces the assertion to the special case where~$G$ is a Lie group.
  The projection $G\to G/H$ is a submersion in this case and hence has
  local smooth sections.  They together with smooth partitions of
  unity yield the assertion.
\end{proof}

The modular function $\modular_G\colon G\to\R\inv_+$ is a continuous
group homomorphism.  We define it by the convention
$\modular_G(x)\,dg=d(gx)$.  We have $\modular_G\in\CINF(G)$ because
group homomorphisms are smooth maps and the identical function
$\R\inv_+\to\R$ is a smooth function on~$\R\inv_+$.  Hence
multiplication by~$\modular_G$ is a bornological isomorphism on
$\CCINF(G,V)$ and $\CINF(G,V)$.

If $H\subseteq G$ is an \emph{open} subgroup, then there are
bornological embeddings
$$
\CCINF(H,V)\to\CCINF(G,V),
\qquad
\CINF(H,V)\to\CINF(G,V),
$$
which extend a function on~$H$ by~$0$ outside~$H$.  Its range is the
space of functions supported in~$H$ and thus a retract.  Let
$(G_i)_{i\in I}$ be a directed family of open subgroups of~$G$ with
$G=\bigcup G_i$.  Then $\CCINF(G,V)$ is the strict inductive limit of
the subspaces $\CCINF(G_i,V)$.

We have
\begin{equation}
  \label{eq:CCINF_product}
  \CCINF(G_1\times G_2) \cong
  \CCINF(G_1)\hot\CCINF(G_2) \cong
  \CCINF(G_1,\CCINF(G_2))
\end{equation}
for all locally compact groups $G_1$ and~$G_2$ because the
corresponding result holds for manifolds and the bornological tensor
product commutes with direct limits.  The spaces $\CINF(G\times G)$,
$\CINF(G)\hot\CINF(G)$ and $\CINF(G,\CINF(G))$ agree if~$G$ is a Lie
group, but not for arbitrary~$G$.

Let $(G_i)_{i\in I}$ be a set of locally compact groups and let
$K_i\subseteq G_i$ be compact open subgroups for all $i\in I\setminus
F_0$ with some \emph{finite} set of exceptions~$F_0$.  For each finite
subset $F\subseteq I$ containing~$F_0$ the direct product
$$
G(F) \defeq \prod_{i\in I\setminus F} K_i \times \prod_{i\in F} G_i
$$
is a locally compact group.  For $F_1\subseteq F_2$ the group
$G(F_1)$ is an open subgroup of $G(F_2)$.  The \emph{restricted direct
  product} $\rprod_{i\in I} (G_i,K_i)$ is the direct union of these
groups.  The characteristic function of $K_i\subseteq G_i$ is a
distinguished element of $\CCINF(G_i)$.  The \emph{(restricted) tensor
  product} of the spaces $\CCINF(G_i)$ with respect to these
distinguished vectors is defined as follows.  For each finite subset
$F\subseteq I$ containing~$F_0$, consider the completed tensor product
$\bigotimes_{i\in F} \CCINF(G_i)$.  We have a map between the
associated tensor products for $F_1\subseteq F_2$ that inserts the
factor $1_{K_i}$ for $i\in F_2\setminus F_1$.  The tensor product is
the direct limit of the resulting (strict) inductive system.  It is
straightforward to show that
\begin{equation}
  \label{eq:CCINF_restricted_product}
  \CCINF \biggl( \rprod_{i\in I} (G_i, K_i) \biggr) \cong
  \bigotimes_{i\in I} (\CCINF(G_i), 1_{K_i}).
\end{equation}

\subsection{Multiplication and convolution}
\label{sec:CCINF_algebra}

The pointwise product of smooth functions and of smooth functions with
distributions is defined in the usual way.  All resulting bilinear
maps are clearly bounded.

The group law of~$G$ gives rise to a comultiplication
$$
\Delta\colon \CINF(G)\to\CINF(G\times G),
\qquad
\Delta f(g,h)\defeq f(gh)
$$
We do not have $\CINF(G\times G)=\CINF(G)\hot\CINF(G)$ in general.
The resulting problem with the convolution of distributions is fixed
by the following lemma:

\begin{lemma}  \label{lem:product_distributions}
  There is a unique bounded bilinear map
  $$
  \CINF'(G)\times\CINF'(G)\to\CINF'(G\times G),
  \qquad
  (D_1,D_2)\mapsto D_1\otimes D_2,
  $$
  such that
  \begin{align*}
    \braket{D_1\otimes D_2}{f_1\otimes f_2} &=
    \braket{D_1}{f_1}\cdot \braket{D_2}{f_2},
    \\
    (f_1\otimes f_2)\cdot(D_1\otimes D_2) &=
    f_1\cdot D_1\otimes f_2\cdot D_2,
  \end{align*}
  for all $D_1,D_2\in\CINF'(G)$, $f_1,f_2\in\CINF(G)$.
  
  There is a unique bounded linear map
  $$
  \CINF'(G/H)\to \Hom(\CINF(G/H,V), V),
  \qquad D\mapsto D_V,
  $$
  such that
  $$
  \braket{D_V}{f\otimes v}= \braket{D}{f}\cdot v,
  \qquad
  f\cdot D_V = (f\cdot D)_V
  $$
  for all $D\in\CINF'(G/H)$, $f\in\CINF(G/H)$, $v\in V$.
\end{lemma}

\begin{proof}
  Fix $D_1,D_2$ with support contained in some compact subset
  $S\subseteq G$.  There exists $\phi\in\CCINF(G)$ with $\phi=1$ in a
  neighborhood of~$S$.  Hence $\phi\cdot D_j=D_j$ for $j=1,2$.
  Therefore, we must put $\braket{D_1\otimes D_2}{f}\defeq
  \braket{D_1\hot D_2}{(\phi\otimes\phi)\cdot f}$.  The right hand
  side is well-defined because $(\phi\otimes\phi)\cdot f$ has compact
  support and $\CCINF(G\times G)\cong\CCINF(G)\hot\CCINF(G)\subseteq
  \CINF(G)\hot\CINF(G)$.  It is straightforward to see that this
  definition does not depend on~$\phi$ and has the required
  properties.
  
  The map~$D_V$ is defined similarly.  There is $\phi\in\CCINF(G/H)$
  with $\phi\cdot D=D$.  We must have $\braket{D_V}{f}\defeq
  D\hot\ID_V(\phi\cdot f)$ for all $f\in\CINF(G/H,V)$.  The right hand
  side is defined because $\phi\cdot
  f\in\CCINF(G/H,V)\cong\CCINF(G/H)\hot V$.
\end{proof}

We define the convolution of two compactly supported distributions by
$$
\braket{D_1\ast D_2}{f} \defeq \braket{D_1\otimes D_2}{\Delta f}
$$
for all $f\in\CINF(G)$.  This turns $\CINF'(G)$ into a bornological
algebra.  A similar trick allows to define the convolution of a
compactly supported distribution with an arbitrary distribution.  All
these bilinear maps are evidently bounded.

Fix a left Haar measure~$dg$ on~$G$.  Then we embed
$\CINF(G)\subseteq\CCINF'(G)$ by the usual map $f\mapsto f\,dg$.  We
define convolutions involving smooth functions in such a way that
$f_1\,dg \ast f_2\,dg = (f_1\ast f_2)\,dg$, $D\ast (f\,dg) = (D\ast
f)\,dg$ and $(f\,dg\ast D)=(f\ast D)\,dg$.  It is straightforward to
verify that this defines bounded bilinear maps taking values in
$\CINF(G)$ provided one factor has compact support, and taking values
in $\CCINF(G)$ if both factors have compact support.  In particular,
$\CCINF(G)$ becomes a bornological algebra and a bimodule over
$\CINF'(G)$.

The antipode $\tilde{f}(g)\defeq f(g^{-1})$ on $\CINF(G)$ gives rise
by transposition to an antipode on $\CINF'(G)$, which is a bounded
anti\brd{}homomorphism with respect to convolution.  Its restriction to the
ideal $\CCINF(G)\subseteq\CINF'(G)$ is given by
\begin{equation}
  \label{eq:CCINF_antipode}
  (\tilde{f}^{(1)})(g)\defeq f(g^{-1})\modular_G(g)^{-1}
\end{equation}
because $d(g^{-1})=\modular_G(g^{-1})\,dg$.  This is a bounded
anti\brd{}homomorphism on $\CCINF(G)$, which we use to turn right
$\CCINF(G)$\brd{}modules into left modules and vice versa.

\section{Smooth representations of locally compact groups}
\label{sec:smooth_representation}

We shall use the following notation and conventions.  Let~$G$ be a
locally compact group and let~$V$ be a (complete convex) bornological
vector space.  The space $\End(V)\defeq\Hom(V,V)$ of bounded linear
operators on~$V$ is a (complete convex) bornological algebra.  Let
$\Aut(V)$ be the multiplicative group of invertible elements in
$\End(V)$.  A \emph{group representation} of~$G$ on~$V$ is a group
homomorphism $\pi\colon G\to\Aut(V)$.  Thus we always assume~$G$ to
act by bounded linear operators.  We write $\pi(g)=\pi_g$ and
$\pi_g(v)=\pi(g,v)=g\cdot v$.  Let $\Map(G,V)\defeq \prod_{g\in G} V$
be the space of all functions from~$G$ to~$V$.  The \emph{adjoint}
of~$\pi$ is the bounded linear map $\pi_\ast \colon V \to \Map(G,V)$
defined by $\pi_\ast(v)(g)\defeq \pi(g,v)$.  We let~$G$ act on
$\Map(G,V)$ by the right regular representation~$\rho$ defined
in~\eqref{eq:biregular_representation}.  Then~$\pi_\ast$ is
$G$\nbd{}equivariant.

\begin{definition}  \label{def:smooth_rep}
  The representation~$\pi$ is called \emph{smooth} if~$\pi_\ast$ is a
  bounded map into $\CINF(G,V)$.
\end{definition}

\subsection{First properties of smooth representations}
\label{sec:smooth_first}

\begin{lemma}  \label{lem:smoothness_local}
  The representation~$\pi$ is smooth if and only if $Wf(x)\defeq
  x\cdot f(x)$ defines an element of $\Aut(\CCINF(G,V))$.  Even more,
  $\pi$ is already smooth if
  $$
  W_\phi\colon
  V \overset{\phi_\ast}\longrightarrow \CCINF(G,V)
  \overset{W}\longrightarrow \Map(G,V),
  \qquad
  v \mapsto [g\mapsto \phi(g)\pi(g,v)],
  $$
  is a bounded linear map into $\CCINF(G,V)$ for some non-zero
  $\phi\in\CCINF(G)$.
\end{lemma}

\begin{proof}
  We have $W_\phi(v)=W(\phi\otimes v)= M_\phi\pi_\ast(v)$,
  where~$M_\phi$ denotes the operator of pointwise multiplication
  by~$\phi$ on $\CCINF(G,V)$.  It follows from the definition of
  $\CINF(G,V)$ that~$\pi$ is smooth if and only if~$W_\phi$ is a
  bounded linear map into $\CCINF(G,V)$ for all~$\phi$.  This is
  equivalent to~$W$ being a bounded linear map.  If~$W$ is bounded, so
  is its inverse $W^{-1}f(x)\defeq x^{-1} f(x)$.  Hence~$W$ belongs to
  $\Aut(\CCINF(G,V))$ if and only if~$\pi$ is smooth.
  
  It remains to prove that~$W_\phi$ is a bounded map into
  $\CCINF(G,V)$ for all $\phi\in\CCINF(G)$ once this happens for a
  single $\phi\neq0$.  Let $X\subseteq\CCINF(G)$ be the subspace of
  all~$\phi$ for which~$W_\phi$ is a bounded map into $\CCINF(G,V)$.
  Clearly, $X$ is an ideal for the pointwise product.  Since $\pi(g)$
  is bounded for all $g\in G$, the operator~$W_\phi$ is bounded if and
  only if $W_{\rho_g\phi}$ is bounded.  Hence for all $g\in G$ there
  is $\phi\in X$ with $\phi(g)\neq0$.  Since~$X$ is an ideal, we get
  $X=\CCINF(G)$.
\end{proof}

\begin{corollary}  \label{cor:smooth_open_subgroup}
  Let $H\subseteq G$ be an open subgroup.  Then a representation
  of~$G$ is smooth if and only if its restriction to~$H$ is smooth.
  Any representation of a discrete group is smooth.
\end{corollary}

\begin{lemma}  \label{lem:CCINF_smooth}
  Let $H\subseteq G$ be a closed subgroup.  The left and right regular
  representations of~$G$ on $\CCINF(G/H,V)$ and $\CCINF(H\backslash
  G,V)$ are smooth.
\end{lemma}

\begin{proof}
  We observed after Lemma~\ref{lem:CCINF_functoriality} that the map
  $G\times G/H\to G\times G/H$ that sends $(x,yH)$ to $(x,xyH)$ is
  smooth.  Since it is also proper, it induces a bounded linear
  operator on $\CCINF(G,\CCINF(G/H,V))\cong\CCINF(G\times G/H,V)$.
  This is the operator~$W$ of Lemma~\ref{lem:smoothness_local} for the
  left regular representation~$\lambda$ on $\CCINF(G/H,V)$.
  Hence~$\lambda$ is smooth.  Similarly, the right regular
  representation on $\CCINF(H\backslash G,V)$ is smooth.
\end{proof}

The regular representations on $\CINF(G,V)$ usually fail to be
smooth.  See Section~\ref{sec:smooth_diff} for some positive
results on $\CINF(G,V)$.

The \emph{integrated form} of a smooth representation~$\pi$ is the
bounded homomorphism
$$
\IN\pi\colon \CINF'(G)\to\End(V),
\qquad
\IN\pi(D)(v)\defeq
D_V\bigl(\pi_\ast(v)\bigr).
$$
The operator $D_V\colon \CINF(G,V)\to V$ is defined in
Lemma~\ref{lem:product_distributions}.  We evidently have
$\IN\pi(\delta_g) =\pi_g$, so that~$\IN\pi$ extends~$\pi$.  We omit
the straightforward proof that~$\IN\pi$ is an algebra homomorphism.
Let $\UniEnvel(G)\subseteq\CINF'(G)$ be the subalgebra of
distributions supported at~$1_G$.  If~$G$ is a Lie group with Lie
algebra~$\LG$, then $\UniEnvel(G)$ is the universal enveloping algebra
of~$\LG$.  Restricting~$\IN\pi$ to $\LG\subseteq\UniEnvel(G)$, we
obtain a Lie algebra representation $D\pi\colon \LG\to\End(V)$.  We
call~$D\pi$ the \emph{differential of~$\pi$}.

\subsection{Permanence properties of smooth representations}
\label{sec:smooth_permanence}

\begin{lemma}  \label{lem:smooth_hereditary}
  Smoothness is hereditary for subrepresentations and quotients,
  direct limits and finite inverse limits (that is, inverse limits of
  finite diagrams).
\end{lemma}

\begin{proof}
  Let $K\into E\prto Q$ be a bornological extension of representations
  of~$G$.  Consider the diagram
  $$
  \xymatrix{
    {K\;} \ar@{>->}[r] \ar@{.>}[d] &
       E \ar@{->>}[r] \ar[d] &
          Q \ar@{.>}[d] \\
    {\CINF(G,K)\;}  \ar@{>->}[r] &
      {\CINF(G,E)}  \ar@{->>}[r] &
        {\CINF(G,Q).} \\
  }
  $$
  The middle vertical map is the adjoint of the representation on~$E$.
  The bottom row is a bornological extension as well by
  Proposition~\ref{pro:CINF_exact}.  Since the composition $K\to
  E\to\CINF(G,E)\to\CINF(G,Q)$ vanishes, the dotted arrows exist.
  They are the adjoints of the induced representations on $K$ and~$Q$.
  Hence $K$ and~$Q$ are smooth representations as well.  It is trivial
  to verify that direct sums of smooth representations are again
  smooth.  Since direct limits are quotients of direct sums and
  inverse limits are subspaces of direct products, we obtain the
  asserted smoothness for direct limits and finite inverse limits.
\end{proof}

\begin{remark}  \label{rem:smooth_not_hereditary}
  Infinite direct products of smooth representations may fail to be
  smooth.  The class of smooth representations is \emph{not} closed
  under extensions.  A simple counterexample is the representation
  of~$\R$ on~$\C^2$ by
  $$
  t\mapsto
  \begin{pmatrix} 1 & \phi(t) \\ 0 & 1
  \end{pmatrix}
  $$
  for some discontinuous group homomorphism $\phi\colon \R\to\R$.
\end{remark}

\begin{lemma}  \label{lem:smoothness_group_functorial}
  Let $\phi\colon H\to G$ be a continuous group homomorphism and let
  $\pi\colon G\to \End(V)$ be a group representation.  If~$\pi$ is a
  smooth representation of~$G$, then $\pi\circ\phi$ is a smooth
  representation of~$H$.  In particular, restrictions of smooth
  representations to closed subgroups remain smooth.  If~$\phi$ is an
  open surjection, then the converse holds.  That is, a representation
  of a quotient group $H/N$ is smooth if and only if it is smooth as a
  representation of~$H$.
\end{lemma}

\begin{proof}
  The smoothness of $\pi\circ\phi$ follows from the functoriality of
  $\CINF(G,V)$ for continuous group homomorphisms.  If~$\phi$ is an
  open surjection, it is isomorphic to a quotient map $\phi\colon H\to
  H/N$.  The map $\phi^\ast\colon \CINF(H/N,V)\to\CINF(H,V)$ is a
  bornological isomorphism onto its range by
  Lemma~\ref{lem:CINF_invariant}.  Hence $\pi\circ\phi$ is smooth if
  and only if~$\pi$ is.
\end{proof}

The \emph{external tensor product} $\pi_1\boxtimes\pi_2$ of two
representations $\pi_j\colon G_j\to\Aut(V_j)$, $j=1,2$, is the tensor
product representation of $G_1\times G_2$ on $V_1\hot V_2$.  If
$G_1=G_2$, the \emph{internal tensor product} $\pi_1\hot\pi_2$ is the
restriction of $\pi_1\boxtimes\pi_2$ to the diagonal $G\subseteq
G\times G$.  Let $(G_i)_{i\in I}$ and $(K_i)_{i\in I\setminus F_0}$ be
the data for a restricted direct product of groups.  Let $\pi_i\colon
G_i\to\Aut(V_i)$ be representations of~$G_i$ and let $\xi_i\in V_i$ be
$K_i$\nbd{}invariant for all but finitely many $i\in I$.  Then we can
form the restricted tensor product $\bigotimes_{i\in I} (V_i,\xi_i)$
and let $\rprod_{i\in I} (G_i,K_i)$ act on it in the evident fashion.
We call this the \emph{restricted (external) tensor product
  representation}.  This recipe is frequently used to construct
representations of adelic groups.

\begin{lemma}  \label{lem:smooth_direct_product}
  A representation of a direct product group is smooth if and only if
  its restrictions to the factors are smooth.  Restricted external
  tensor products and external and internal tensor products of smooth
  representations remain smooth.
\end{lemma}

\begin{proof}
  The straightforward proof of the first assertion is left to the
  reader.  Consider a restricted direct product $G=\rprod (G_i,K_i)$
  and a restricted tensor product representation $\bigotimes_{i\in I}
  (V_i,\xi_i)$ as above.  We have
  $$
  \CCINF\bigl(G,\bigotimes (V_i,\xi_i)\bigr) \cong
  \bigotimes (\CCINF(G_i),1_{K_i}) \hot
  \bigotimes (V_i,\xi_i) \cong
  \bigotimes (\CCINF(G_i,V_i),1_{K_i}\otimes\xi_i).
  $$
  The restricted tensor product is functorial for families of maps
  $V_i\to V_i$ preserving the distinguished vectors.  Since the
  operator~$W$ of Lemma~\ref{lem:smoothness_local} is induced from the
  analogous operators for the factors, we get the assertion for
  restricted direct products.  This implies the smoothness of finite
  external tensor products and hence also of internal tensor products
  by Lemma~\ref{lem:smoothness_group_functorial}.
\end{proof}

\subsection{Some constructions with representations}
\label{sec:constructions_rep}

\begin{definition}  \label{def:smoothen_rep}
  The \emph{smoothening} of a representation $\pi\colon G\to\Aut(V)$
  is
  $$
  \smooth_G V\defeq
  \{ f\in\CINF(G,V)\mid
  \text{$f(g)= g\cdot f(1)$ for all $g\in G$}\},
  $$
  equipped with the subspace bornology, the right regular
  representation and the map $\iota_V\colon \smooth_G V\to V$ defined
  by $\iota_V(f)=f(1)$.
\end{definition}

We frequently drop~$G$ and just write $\smooth(V)$ for the
smoothening.  We write $\smooth_G (V,\pi)$ if it is important to
remember the representation~$\pi$.  A function $f\in\CINF(G,V)$
belongs to $\smooth(V)$ if and only if $f=\pi_\ast\bigl(f(1)\bigr)$.
Therefore, the map~$\iota_V$ is injective and $\smooth(V)$ is
invariant under the right regular representation.  The map~$\iota_V$
is bounded and $G$\nbd{}equivariant.

Let $L\subseteq G$ be a compact neighborhood of the identity.  Recall
that $\CINF(L,V)$ is defined as a quotient of $\CINF(G,V)$ in
Definition~\ref{def:CINF}.  However, since~$L$ is compact, it is also
a quotient of $\CCINF(G,V)$.  Therefore, $\CINF(L,V)\cong\CINF(L)\hot
V$.

\begin{lemma}  \label{lem:smoothen_local}
  The projection $(v,f)\mapsto f|_L$ is a bornological isomorphism
  from $\smooth(V)$ onto the space
  $$
  \smooth_L V\defeq
  \{f\in\CINF(L,V)\mid
  \text{$f(g)=g\cdot f(1)$ for all $g\in L$}\}.
  $$
  In particular, $\smooth_H V\cong\smooth_G V$ if $H\subseteq G$
  is an open subgroup.
\end{lemma}

\begin{proof}
  Restriction to~$L$ is a bounded linear map $p\colon
  \smooth(V)\to\smooth_L V$.  Define $jf(g)\defeq g\cdot f(1)$ for
  all $g\in G$, $f\in\smooth_L V$.  This is a bounded linear map from
  $\smooth_L V$ to $\smooth(V)$ because $j(f)|_{gL}=\pi_g(f)$ and the
  interiors of the sets~$gL$ with $g\in G$ cover~$G$.  Clearly, the
  maps $j$ and~$p$ are inverse to each other.
\end{proof}

\begin{proposition}  \label{pro:smoothen}
  The smoothening of~$V$ is a smooth representation of~$G$.
  If~$W$ is any smooth representation of~$G$, then there is a natural
  isomorphism
  $$
  (\iota_V)_\ast\colon
  \Hom_G(W,V)\cong \Hom_G(W,\smooth(V)).
  $$
\end{proposition}

\begin{proof}
  The map $(\iota_V)_\ast$ is injective because~$\iota_V$ is.  A map
  $T\colon W\to V$ induces a map
  $\CINF(G,T)\colon\CINF(G,W)\to\CINF(G,V)$.  We have
  $\iota_V\circ\CINF(G,T)\circ\pi_\ast^W=T$ and
  $\CINF(G,T)\circ\pi_\ast^W$ maps~$W$ into $\smooth(V)$ if~$T$ is
  equivariant.  Hence $(\iota_V)_\ast$ is also surjective.
  
  It remains to prove the smoothness of $\smooth(V)$.  This requires
  work because the regular representation on $\CINF(G,V)$ may fail to
  be smooth.  Let $L\subseteq G$ be a compact symmetric neighborhood
  of~$1$ and let $L^2\defeq L\cdot L$.  There is a bounded linear map
  $$
  \rho^\ast\colon \CINF(G,V)\to\CINF(G\times G,V),
  \qquad
  \rho^\ast f(g,h)\defeq f(gh).
  $$
  It descends to a bounded map $\CINF(L^2,V)\to\CINF(L\times
  L,V)\cong\CINF(L,\CINF(L,V))$, which maps $\smooth_{L^2}(V)$ into
  $\CINF(L,\smooth_L V)$.  The isomorphism $\CINF(L\times
  L,V)\cong\CINF(L,\CINF(L,V))$ follows immediately from
  $\CINF(L,V)\cong \CINF(L)\hot V$, but it holds only if~$L$ is
  compact.  Using Lemma~\ref{lem:smoothen_local}, we get a bounded map
  $$
  \rho^\ast\colon \smooth(V)\to\CINF\bigl(L,\smooth(V)\bigr),
  \qquad
  \rho^\ast(f)(g)\defeq \rho_g(f),
  $$
  Since~$L$ is a neighborhood of the identity, the smoothness of
  $\smooth(V)$ now follows from Lemma~\ref{lem:smoothness_local}.
\end{proof}

Let~$\NRepG$ be the category of representations of~$G$ on bornological
vector spaces with $G$\nbd{}equivariant bounded linear maps as
morphisms.  Let~$\RepG$ be the full subcategory of smooth
representations.  Proposition~\ref{pro:smoothen} asserts that
$\smooth\colon \NRepG\to\RepG$ is right adjoint to the embedding
$\RepG\subseteq\NRepG$.

Let $H\subseteq G$ be a closed subgroup.  We have an evident
restriction functor $\Res_G^H\colon \NRepG\to\NRepH$, which maps
$\RepG$ into~$\RepH$.  The \emph{smooth induction functor}
$\Ind_H^G\colon \RepH\to\RepG$ is defined as the right adjoint of the
restriction functor.  The following construction shows that it exists.

First we construct a right adjoint to $\Res_G^H\colon
\NRepG\to\NRepH$.  Let
$$
I(V)\defeq
\{v\in \Map(G,V)\mid
  \text{$f(hg)= h\cdot f(g)$ for all $h\in H$, $g\in G$}
\},
$$
equipped with the subspace bornology from $\Map(G,V)$ and the right
regular representation.  A morphism $f\colon \Res_G^H(W)\to V$ in
$\NRepH$ induces a morphism $f_\ast\colon W\to I(V)$ in $\NRepG$ by
$f_\ast(w)(g)\defeq f(gw)$.  Any morphism $W\to I(V)$ is of this form
for a unique morphism~$f$.  That is, $I$ is right adjoint to the
restriction functor $\NRepG\to\NRepH$.  It follows easily that the
functor
$$
\Ind_H^G\colon \RepH\to\RepG,
\qquad V\mapsto \smooth_G I(V),
$$
is right adjoint to the restriction functor $\RepG\to\RepH$.  Any
$G$\nbd{}equivariant map $W\to \Map(G,V)$ for a smooth
representation~$W$ already takes values in $\CINF(G,V)$.  Hence we can
use $\CINF(G,V)$ instead of $\Map(G,V)$ to define of $\Ind_H^G(V)$.
However, we still have to smoothen afterwards because $\CINF(G,V)$ may
fail to be smooth.

The support of a function in $I(V)$ is left-$H$\brd{}invariant and can be
viewed as a subset of $H\backslash G$.  We let $I_c(V)$ be the
subspace of compactly supported functions in $I(V)$, equipped with the
inductive limit bornology over the compact subsets of $H\backslash G$.
We define the \emph{compact induction} functor as
$$
\cInd_H^G\colon \RepH\to\RepG,
\qquad
V\mapsto \smooth_G I_c(V).
$$

\begin{proposition}  \label{pro:cInd_concrete}
  The representation $\cInd_H^G(V)$ is isomorphic to the right regular
  representation of~$G$ on
  $$
  W\defeq \varinjlim {} \{f\in \CINF_0(H\cdot S,V)\mid
  \text{$f(hg)= h\cdot f(g)$ for all $h\in H$, $g\in G$} \},
  $$
  where~$S$ runs through the compact subsets of $H\backslash G$.

  The functor $\cInd_H^G$ preserves direct limits, injectivity of
  morphisms, bornological extensions, linearly split extensions and
  locally linearly split extensions.
\end{proposition}

\begin{proof}
  It is clear that~$W$ is a subrepresentation of $I_c(V)$.
  Furthermore, any map $X\to I_c(V)$ from a smooth representation to
  $I_c(V)$ must factor through~$W$.  We must prove that~$W$ is a
  smooth representation of~$G$.  We do this by realizing it naturally
  as a linearly split quotient of the left regular representation on
  $\CCINF(G,V)$.  Thus the functor $\cInd_H^G$ is a retract of the
  functor $\CCINF(G,\blank)$ if we forget the group representation.
  Hence it inherits its functorial properties listed in
  Proposition~\ref{pro:CCINF_exact}.

  Consider the maps
  \begin{equation}  \label{eq:cInd_concrete}
  \begin{alignedat}{2}
    P&\colon \CCINF(G,V)\to W,
    &\qquad
    Pf(g) &\defeq \int_H h\cdot f(g^{-1}h)\,d_Hh,
    \\
    J&\colon W\to\CCINF(G,V),
    &\qquad
    J f(g) &\defeq f(g^{-1})\cdot\phi(g).
  \end{alignedat}
  \end{equation}
  The map~$P$ is bounded and $G$\nbd{}equivariant.  The map~$J$ is a
  bounded linear left section for~$P$ provided $\supp \phi\cap S\cdot
  H$ is compact for all $S\subseteq G/H$ compact and $\int_H
  \phi(gh)\,d_Hh=1$ for all $g\in G$.  Such a function~$\phi$ clearly
  exists.  As a quotient of the left regular representation on
  $\CCINF(G,V)$, the representation~$W$ is smooth.
\end{proof}

Proposition~\ref{pro:cInd_concrete} easily implies that
\begin{align}
  \label{eq:cInd_regular}
  \cInd_H^G(\CCINF(H,V)) &\cong
  \CCINF(G,V),
  \\
  \label{eq:cInd_trivial}
  \cInd_H^G(\C(1)) &\cong \CCINF(G/H),
\end{align}
where $\C(1)$ denotes the trivial representation of~$H$ on~$\C$ and
all function spaces carry the left regular representation.

It is customary to twist the functors $\Ind_H^G$ and $\cInd_H^G$ by a
modular factor.  Let $\modular_G$ and~$\modular_H$ be the modular
functions of $G$ and~$H$, respectively.  We call the quasi\brd{}character
$\modular_{G:H}\defeq \modular_G\modular_H^{-1}\colon H\to\R\inv_+$
the \emph{relative modular function} of $H\subseteq G$.  For a
representation $\pi\colon H\to\Aut(W)$ of~$H$ and $\alpha\in\R$, we
form the representation $\modular_{G:H}^\alpha\cdot\pi$ on~$W$ and plug
it into $\Ind_H^G$ and $\cInd_H^G$ instead of~$W$ itself.  We call the
resulting functors the \emph{twisted induction and compact induction
  functors}.  The case $\alpha=1/2$ is important because it preserves
unitary representations.

\subsection{Explicit criteria for smoothness}
\label{sec:explicit_smooth_rep}

Let $U\subseteq G$ be an open subgroup which is a projective limit of
Lie groups.  Let~$I$ be a fundamental system of smooth compact
subgroups in~$U$.  For a subgroup $L\subseteq G$ we let
$$
V^L\defeq\{v\in V\mid \text{$gv=v$ for all $g\in L$}\}.
$$
This is a closed linear subspace of~$V$.  The subspaces~$V^k$ for
$k\in I$ form a strict inductive system.  We have $V=\varinjlim V^k$
if and only if any bounded subset of~$V$ is contained in~$V^k$ for
some $k\in I$.

\begin{theorem}  \label{the:explicit_criterion_smooth_I}
  A representation $\pi\colon G\to\Aut(V)$ is smooth if and only if
  $V=\varinjlim_{k\in I} V^k$ and the representation of $U/k$ on~$V^k$
  is smooth for all $k\in I$.
\end{theorem}

\begin{proof}
  Since~$\pi$ is smooth if and only if its restriction to~$U$ is
  smooth we may assume without loss of generality that $G=U$.  We may
  also assume that there be $k_0\in I$ with $k\subseteq k_0$ for all
  $k\in I$.  Fix $\phi\in\CCINF(G/k_0)$ with $\phi(1)\neq0$.  The
  representation~$\pi$ is smooth if and only if the operator~$W_\phi$
  in Lemma~\ref{lem:smoothness_local} is a bounded map from~$V$ to
  $\CCINF(G,V)\cong \varinjlim \CCINF(G/k,V)$.  Evidently, $W_\phi(v)$
  is $k$\nbd{}invariant if and only if $v\in V^k$.  As a result, we
  must have $V=\varinjlim V^k$ if~$\pi$ is smooth.  Suppose now that
  $V=\varinjlim V^k$.  Since smoothness is hereditary for inductive
  limits and subrepresentations, $V$ is smooth if and only if~$V^k$ is
  smooth for all $k\in I$.  Moreover, the representation of~$G$
  on~$V^k$ is smooth if and only if the induced representation of
  $G/k$ is smooth.  This yields the assertion.
\end{proof}

If~$G$ is totally disconnected, the quotients $U/k$ are discrete, so
that any representation of~$U/k$ is smooth.  Therefore, $\pi$ is
smooth if and only if $V=\varinjlim V^k$.  If~$V$ carries the fine
bornology, then the latter holds if and only if each $v\in V$ is
stabilized by some open subgroup.  For arbitrary~$G$ the quotients
$U/k$ are Lie groups.  Hence it remains to describe smooth Lie group
representations.

\begin{theorem}  \label{the:explicit_criterion_smooth_II}
  Let~$G$ be a Lie group and let~$\LG$ be its Lie algebra.  A
  representation $\pi\colon G\to\Aut(V)$ is smooth if and only if it
  satisfies the following conditions:
  \begin{enumerate}[(i)]
  \item the representation is locally equibounded, that is,
    $\pi(K)\subseteq \End(V)$ is equibounded for any compact subset
    $K\subseteq G$;
    
  \item the limits $D\pi(X)(v) \defeq \lim_{t\to0}
    t^{-1}(\exp(tX)\cdot v-v)$ exist for all $v\in V$ and the
    convergence is uniform on bounded subsets of~$V$;
    
  \item for any bounded subset $S\subseteq V$ there is a bounded disk
    $T\subseteq V$ such that $D\pi(X_1)\circ\dots\circ D\pi(X_n)(S)$
    is bounded in~$V_T$ for all $n\in\N$, $X_1,\dots,X_n\in\LG$.
  \end{enumerate}
\end{theorem}

\begin{proof}
  First we show that smooth representations satisfy (i)--(iii).
  Conditions (i) and~(ii) are obvious with $D\pi(X)=\IN\pi(X)$ for
  all $X\in\LG\subseteq\CINF'(G)$.  Let $S\subseteq V$ be bounded and
  let $\phi\in\CCINF(G)$ be such that $\phi=1$ in a neighborhood of
  the identity.  Define $W_\phi(v)(g)\defeq \phi(g)\pi(g,v)$ as in
  Lemma~\ref{lem:smoothness_local}.  The set $W_\phi(S)$ is bounded in
  $\CCINF(G,V)$ and hence in $\CCINF(G,V_T)$ for some bounded disk
  $T\subseteq V$.  This yields~(iii).
  
  Conversely, suppose (i)--(iii) to hold.  We claim that~$\pi$ is
  smooth.  Let $S\subseteq V$ be a bounded complete disk and
  $K\subseteq G$ compact.  Condition~(i) allows us to choose a bounded
  complete disk $S'\subseteq V$ containing $\pi(K)(S)$.  Let
  $S''\subseteq V$ be a bounded complete disk such that the
  convergence in~(ii) is uniform in~$V_{S''}$ for all $v\in S'$.  Such
  a set exists by the definition of uniform convergence.
  Condition~(iii) asserts that there is a bounded complete disk~$T$
  such that $D\pi(X_1)\circ\dots\circ D\pi(X_n)(S'')$ is bounded
  in~$V_T$ for all $n\in\N$, $X_1,\dots,X_n\in\LG$.
  
  We claim that the map $v\mapsto \pi_\ast(v)|_K$ is a bounded linear
  map from~$V_S$ to $\CINF(K,V_T)$.  This claim implies that~$\pi$ is
  smooth.  Since $V_S$ and~$V_T$ are Banach spaces, the claim is
  equivalent to the smoothness of the Banach space valued map
  $\pi\colon K\to\Hom(V_S,V_T)$.  This is what we are going to show.
  The construction of the sets $S',S'',T$ yields the following.  The
  family of operators $\pi(g)\colon V_S\to V_{S'}$ is uniformly
  bounded for $g\in K$.  Let $X_1,\dots,X_n,X\in\LG$.  The operators
  $(\pi(\exp(hX))-\ID)/h\colon V_{S'}\to V_{S''}$ converge
  towards~$D\pi$ in operator norm for $h\to0$.  The operator $A\defeq
  D\pi(X_1)\circ\dots\circ D\pi(X_n)\colon V_{S''}\to V_{T}$ is
  bounded.  Hence
  \begin{displaymath}
    \lim_{h\to 0} A\circ(\pi(\exp(hX)g)-\pi(g))/h =
    A\circ D\pi(X)\circ \pi(g)
  \end{displaymath}
  converges in $\Hom(V_S,V_T)$ and is of the same form as the operator
  $A\circ\pi(g)$.  This means that we can differentiate~$\pi$ with
  respect to right invariant differential operators.  Therefore, $\pi$
  is a $\CONT^\infty$\brd{}map from~$K$ to $\Hom(V_S,V_T)$ as claimed.
\end{proof}

\subsection{Smooth versus differentiable representations}
\label{sec:smooth_diff}

Let~$G$ be a Lie group.  Using the spaces $\CCONT^k(G,V)$ defined in
Section~\ref{sec:preliminaries} instead of $\CCINF(G,V)$, we define
the space $\CONT^k(G,V)$ of $\CONT^k$\brd{}functions $G\to V$ for
$k\in\N\cup\{\infty\}$ as in Definition~\ref{def:CINF}.  We call~$\pi$
a \emph{$\CONT^k$\brd{}representation} if~$\pi_\ast$ is a bounded map
from~$V$ to $\CONT^k(G,V)$.  For $k=0$ and $k=\infty$ we get
\emph{continuous} and \emph{differentiable} representations,
respectively.

\begin{theorem}  \label{the:explicit_criterion_differentiable}
  Let $\pi\colon G\to\Aut(V)$ be a representation of a Lie group~$G$.
  Let~$\LG$ be the Lie algebra of~$G$.  The following statements are
  equivalent:
  \begin{enumerate}[(1)]
  \item the representation~$\pi$ is differentiable;

  \item the representation~$\pi$ is $\CONT^1$;

  \item there is a bounded homomorphism $\IN\pi\colon
    \CINF'(G)\to\End(V)$ extending~$\pi$;

  \item the following two conditions hold:
    \begin{enumerate}[(i)]
    \item the representation is locally equibounded, that is, for all
      compact subsets $K\subseteq G$ the set $\pi(K)\subseteq
      \End(V)$ is equibounded;
    
    \item the limits $D\pi(X)(v) \defeq \lim_{t\to0}
      t^{-1}(\exp(tX)\cdot v-v)$ exist for all $v\in V$ and the
      convergence is  uniform on bounded subsets of~$V$.
    \end{enumerate}
  \end{enumerate}
\end{theorem}

\begin{proof}
  It is clear that (1) implies~(2).  The dual of $\CONT^1(G)$ is a
  subspace of $\CINF'(G)$.  It generates $\CINF'(G)$ as a bornological
  algebra in the sense that any bounded subset of $\CINF'(G)$ is
  contained in $S^n$ for a bounded subset $S\subseteq \CONT^1(G)'$.  A
  $\CONT^1$\brd{}representation gives rise to a bounded linear map
  $\CONT^1(G)'\to\End(V)$, which we can then extend to an algebra
  homomorphism on all of $\CINF'(G)$.  Hence~(2) implies~(3).  The set
  of~$\delta_g$, $g\in K$, is bounded in $\CINF'(G)$ and we have
  convergence $t^{-1}(\delta_{\exp(tX)}-\delta_1)\to X$ in $\CINF'(G)$
  for all $X\in\LG$.  Hence~(3) implies~(4).  The proof of the
  implication (4)$\Longrightarrow$(1) is similar to the proof of
  Theorem~\ref{the:explicit_criterion_smooth_II} and therefore
  omitted.
\end{proof}

Conditions (i) and~(ii) above are the same as in
Theorem~\ref{the:explicit_criterion_smooth_II}.  Thus the only
difference between smoothness and differentiability is condition~(iii)
of Theorem~\ref{the:explicit_criterion_smooth_II}.

\begin{remark}  \label{rem:distributions_smooth}
  It follows immediately from
  Theorem~\ref{the:explicit_criterion_differentiable} that the regular
  representations on $\CINF'(G)$ and $\CCINF'(G)$ are differentiable.
  However, these representations are not smooth.  One can verify
  directly that the third condition of
  Theorem~\ref{the:explicit_criterion_smooth_II} fails.  It is also
  clear that they are not essential as modules over $\CCINF(G)$
  because the convolution of a smooth function with a distribution is
  already a smooth function.
\end{remark}

\begin{proposition}  \label{pro:CINF_smoothen}
  Let~$G$ be a locally compact group that is countable at infinity and
  let~$V$ be a metrizable bornological vector space.  Let~$I$ be a
  fundamental system of smooth compact subgroups in~$G$.  Then
  \begin{align*}
    \smooth_G (\CINF(G,V),\lambda) &\cong
    \varinjlim_{k\in I} \CINF(k\backslash G,V) \cong
    \varinjlim_{k\in I} \CINF(k\backslash G)\hot V;
    \\
    \smooth_G (\CINF(G,V),\rho) &\cong
    \varinjlim_{k\in I} \CINF(G/k,V) \cong
    \varinjlim_{k\in I} \CINF(G/k)\hot V;
    \\
    \smooth_{G\times G} (\CINF(G,V),\lambda\boxtimes\rho) &\cong
    \varinjlim_{k\in I} \CINF(G//k,V) \cong
    \varinjlim_{k\in I} \CINF(G//k)\hot V.
  \end{align*}
\end{proposition}

\begin{proof}
  We only compute the smoothening of the left regular representation,
  the other cases are similar.  Let $U\subseteq G$ be an open almost
  connected subgroup.  We can assume all $k\in I$ to be normal
  subgroups of~$U$.  Let $k\in I$.  Since $V$ is metrizable and
  $k\backslash G$ is countable at infinity,
  Lemma~\ref{lem:CINF_metrizable} yields $\CINF(k\backslash
  G,V)\cong\CINF(k\backslash G)\hot V$ and hence the last isomorphism.
  The space $\CINF(k\backslash G)\hot V$ is metrizable as well.  Hence
  there is no difference between smooth and differentiable Lie group
  representations on this space by
  Proposition~\ref{pro:smooth_metrizable}.  Since $\CINF'(U/k)$
  evidently acts on $\CINF(k\backslash G)\hot V$ by convolution, we
  conclude that $U/k$ acts smoothly on $\CINF(k\backslash G)\hot V$
  for all $k\in I$.  Therefore, $X\defeq \varinjlim \CINF(k\backslash
  G) \hot V$ is a smooth representation of~$G$ by
  Theorem~\ref{the:explicit_criterion_smooth_I}.  Since $W=\varinjlim
  W^k$ for any smooth representation, it is clear that any bounded
  $G$\nbd{}equivariant map $W\to\CINF(G,V)$ factors through~$X$.
  Hence~$X$ is the smoothening of $\CINF(G,V)$.
\end{proof}

The assertion of the proposition becomes false if~$G$ fails to be
countable at infinity or if~$V$ fails to be metrizable.

\begin{proposition}  \label{pro:smooth_metrizable}
  Differentiable Lie group representations on metrizable bornological
  vector spaces are smooth.
\end{proposition}

\begin{proof}
  This follows immediately from
  Lemma~\ref{lem:differentiable_metrizable}.
\end{proof}

\subsection{Smooth representations on topological vector spaces}
\label{sec:smooth_topological}

Let~$G$ be a Lie group and let~$V$ be a complete locally convex
topological vector space.  Let $\End(V)$ be the algebra of continuous
linear operators on~$V$ and let $\Aut(V)$ be its multiplicative
group.  We equip $\End(V)$ with the equicontinuous bornology, so that
it becomes a bornological algebra.  There is a topological analogue of
the space $\CINF(G,V)$.  A representation $\pi\colon G\to\Aut(V)$ is
called \emph{smooth} if its adjoint is a continuous linear map
$\pi_\ast\colon V\to\CINF(G,V)$
(see~\cite{Bruhat:Representations_induites}).  The following criterion
is similar to the criterion for differentiable representations in
Theorem~\ref{the:explicit_criterion_differentiable}.

\begin{proposition}  \label{pro:smooth_topological}
  The representation~$\pi$ is smooth if and only if it can be extended
  to a bounded homomorphism $\IN\pi\colon \CINF'(G)\to\End(V)$.
\end{proposition}

\begin{proof}
  First suppose~$\pi$ to be smooth.  We let $D\in\CINF'(G)$ act on~$V$
  as usual by $\IN\pi(D)(v)\defeq \braket{D\prot\ID}{\pi_\ast(v)}$.
  This is defined because $\CINF(G,V)\cong\CINF(G)\prot V$ is
  Grothendieck's projective tensor product
  (\cite{Grothendieck:Produits_Tensoriels}).  Let $S\subseteq
  \CINF'(G)$ be bounded.  Then~$S$ is an equicontinuous set of linear
  functionals on $\CINF(G)$ because $\CINF(G)$ is a Fréchet space.
  Hence $\IN\pi(S)$ is equicontinuous as well.  Suppose conversely
  that $\IN\pi\colon \CINF'(G)\to\End(V)$ is a bounded homomorphism
  extending~$\pi$.  Then the family of operators~$\pi_g$ for~$g$ in a
  compact subset of~$G$ is equicontinuous and $t^{-1}(\exp(tX)\cdot
  v-v)\to \IN\pi(X)(v)$ in the strong operator topology for $t\to0$.
  This implies that~$\pi$ is smooth,
  see~\cite{Bruhat:Representations_induites}.
\end{proof}

We now equip~$V$ with the von Neumann bornology, which consists of the
subsets of~$V$ that are absorbed by each neighborhood of zero.  Any
equicontinuous family of operators on~$V$ is equibounded.  Hence a
topologically smooth representation is bornologically differentiable.
The converse implication holds if~$V$ is ``bornological'', that is, a
subset that absorbs all von Neumann bounded subsets is already a
neighborhood of zero.  In that case an equibounded set of linear maps
is equicontinuous as well.  Thus topologically smooth representations
on bornological topological vector spaces are the same as
bornologically differentiable representations with respect to the von
Neumann bornology.

Next we consider the precompact bornology.  Let $\PREC(V)$ be~$V$
equipped with the precompact bornology.  Let~$\pi$ be topologically
smooth.  Since any bounded subset of $\CINF'(G)$ is bornologically
compact, the set of operators $\IN\pi(S)$ for bounded
$S\subseteq\CINF'(G)$ is even bornologically relatively compact for
the equicontinuous bornology on $\End(V)$.  This implies that
$\IN\pi(S)(T)$ is again precompact for precompact~$T$, that is,
$\IN\pi$ is bounded for the equibounded bornology on $\End(\PREC(V))$.
The converse implication holds if a subset of~$V$ that absorbs all
precompact subsets is already a neighborhood of zero.  For instance,
this is the case if~$V$ is a Fréchet space.

As a result, the topological notion of smooth representation is
equivalent to the bornological notion of differentiable representation
under mild hypotheses on the topology of~$V$.  However,
condition~(iii) of Theorem~\ref{the:explicit_criterion_smooth_II} will
usually be violated.

Analogous assertions for continuous representations are false
unless~$V$ is a Fréchet space.  For instance, if~$V$ is a continuous
representation on a Banach space, then the induced representation on
the dual space~$V'$ is weakly continuous but usually not norm
continuous.  However, the weak and the norm topology on~$V'$ have the
same von Neumann bornology.

\begin{theorem}  \label{the:Lie_group_on_Frechet_smooth}
  Let $\pi\colon G\to\Aut(V)$ be a group representation of a Lie group
  on a Fréchet space.  Then the following are equivalent:
  \begin{enumerate}
  \item $\pi$ is smooth as a representation on a topological vector
    space;
  \item $\pi$ is smooth with respect to the von Neumann bornology;
  \item $\pi$ is smooth with respect to the precompact bornology.
  \end{enumerate}
\end{theorem}

\begin{proof}
  A subset of~$V$ that absorbs all null sequences is already a
  neighborhood of zero.  Hence the above discussion shows that
  topological smoothness is equivalent to bornological
  differentiability for either the von Neumann or the precompact
  bornology.  Since both bornologies on~$V$ are metrizable, the
  assertion now follows from Proposition~\ref{pro:smooth_metrizable}.
\end{proof}

\begin{proposition}  \label{pro:Lie_group_on_Frechet_smoothening}
  Let~$V$ be a Fréchet space equipped with the precompact or von
  Neumann bornology and let~$G$ be a Lie group.  Let $\pi\colon
  G\to\Aut(V)$ be a representation.  Then the smoothening of~$V$ is a
  Fréchet space with the precompact or the von Neumann bornology,
  respectively.  If~$V$ is nuclear, so is $\smooth(V)$.
\end{proposition}

\begin{proof}
  Let~$W$ be the Fréchet space of smooth functions $G\to V$ in the
  usual topological sense, equipped with the precompact or von Neumann
  bornology, respectively.  It is shown in~\cite{Meyer:Born_Top} that
  $\CINF(G,V)=W$ as bornological vector spaces, for both bornologies.
  Here we use that the bornologies of locally uniform boundedness and
  locally uniform continuity on $\CINF(G,V)$ coincide.  Since
  $\smooth(V)$ is a closed subspace of $\CINF(G,V)\cong W$, it is a
  Fréchet space as well.  Furthermore, if~$V$ is nuclear, so is~$W$
  and hence its subspace $\smooth(V)$.
\end{proof}

\section{Essential modules versus smooth representations}
\label{sec:essential_modules}

Let~$G$ be a locally compact group.  We are going to identify the
category of smooth representations of~$G$ with the category of
essential modules over the convolution algebra $\CCINF(G)$.  First we
introduce the appropriate notion of an approximate identity in a
bornological algebra and define the notion of an essential module.
Then we compare essential modules over $\CCINF(G)$ with smooth
representations of~$G$.  Finally, we investigate analogues of the
smoothening, restriction, compact induction and induction functors
for representations.

\subsection{Approximate identities and essential modules}
\label{sec:approxid}

\begin{definition}  \label{def:approximate_identity}
  Let~$A$ be a bornological algebra.  We say that~$A$ has an
  \emph{approximate identity} if for each bornologically compact
  subset $S\subseteq A$ there is a sequence $(u_n)_{n\in\N}$ in~$A$
  such that $u_n\cdot x$ and $x\cdot u_n$ converge to~$x$
  uniformly for $x\in S$.
\end{definition}

A subset of a bornological vector space~$V$ is \emph{bornologically
compact} if it is a compact subset of~$V_T$ for some bounded complete
disk $T\subseteq V$.  The uniform convergence in the above definition
means that there is a bounded complete disk $T\subseteq A$ such that
$u_nx$ and $xu_n$ converge to~$x$ uniformly for $x\in S$ in the Banach
space~$V_T$.

Since we may take a different sequence $(u_n)$ for each bornologically
compact subset, we are really considering a net $(u_{n,S})$ in~$A$,
indexed by pairs $(S,n)$ where $S\subseteq A$ is bornologically
compact and $n\in\N$.  It is more convenient to work with sequences as
in Definition~\ref{def:approximate_identity}, however.  The above
definition is related to the usual notion of an approximate identity
in a Banach algebra:

\begin{lemma}  \label{lem:approxid_Banach}
  Let~$A$ be a Banach algebra with a (multiplier) bounded approximate
  identity in the usual sense.  Then~$A$ equipped with the von
  Neumann or precompact bornology has an approximate identity in the
  sense of Definition~\ref{def:approximate_identity}.
\end{lemma}

\begin{proposition}  \label{pro:CCINF_approximate_identity}
  The bornological algebra $\CCINF(G)$ has an approximate identity for
  any locally compact topological group~$G$.
\end{proposition}

\begin{proof}
  Let $U\subseteq G$ be open and almost connected.  Any element of
  $\CCINF(G)$ can be written as a finite sum of elements of the form
  $\delta_g\ast f$ or of elements of the form $f\ast\delta_g$ with
  $g\in G$, $f\in\CCINF(U)$.  Therefore, it suffices to construct an
  approximate identity for $\CCINF(U)$.  Let~$I$ be a fundamental
  system of smooth compact subgroups of~$U$.  Since
  $\CCINF(U)=\varinjlim \CCINF(U/k)$, it suffices to construct
  approximate identities in $\CCINF(U/k)$.  Consequently, we may
  assume~$G$ to be a Lie group.
  
  Let $(u_n)_{n\in\N}$ be a sequence in $\CCINF(G)$ with
  $$
  \lim_{n\to\infty} \int_G u_n(g)\,dg=1,
  \qquad
  \lim_{n\to\infty} \supp u_n=\{1\}.
  $$
  The latter condition means that the support of~$u_n$ is
  eventually contained in any neighborhood of~$1$.  We claim that
  $(u_n)$ is an approximate identity for any bounded subset
  $S\subseteq\CCINF(G)$.  We only check the convergence $u_n\ast f\to
  f$.  The convergence $f\ast u_n\to f$ is proved similarly, using
  that $\lim \int_G u_n(g^{-1})\,dg=1$ as well.
  
  There is a compact subset $K\subseteq G$ such that~$f$ and $f\ast
  u_n$ are supported in~$K$ for all $f\in S$, $n\in\N$.  Hence we are
  working in the nuclear Fréchet space $\CINF_0(K)$.  It is
  straightforward to see that $u_n\ast f$ converges to~$f$ with
  respect to the topology of $\CINF_0(K)$, even uniformly for $f\in
  S$.  Since $\CINF_0(K)$ is a Fréchet space equipped with the von
  Neumann bornology, the topological and bornological notions of
  uniform convergence of a sequence of operators on precompact subsets
  in $\CINF_0(K)$ are equivalent (see~\cite{Meyer:Born_Top}).
  Hence~$(u_n)$ is a left approximate identity in the sense of
  Definition~\ref{def:approximate_identity}.
\end{proof}

Let~$V$ be a right and~$W$ a left bornological $A$\nbd{}module.  Then we
define $V\hot_A W$ as the cokernel of the map
$$
b'_1\colon V\hot A\hot W\to V\hot W,
\qquad
v\otimes a\otimes w\mapsto va\otimes w-v\otimes aw.
$$
That is, we divide $V\hot W$ by the \emph{closure} of the range
of~$b'_1$.  For $V=A$ we also consider the map $b'_0\colon A\hot W\to
W$, $a\otimes w\mapsto aw$.  Since $b'_0\circ b'_1=0$, the map~$b'_0$
descends to a map $A\hot_A W\to W$.  If~$V$ is a $B$\nbd{}$A$\brd{}bimodule
and~$W$ a left $A$\nbd{}module, then $V\hot_A W$ is a left $B$\nbd{}module
in an obvious fashion.  In particular, $A\hot_A W$ is a left
$A$\nbd{}module and the map $A\hot_A W\to W$ is a module homomorphism.

\begin{lemma}  \label{lem:approxid_hota}
  Let~$A$ be a bornological algebra with an approximate identity and
  let~$W$ be a bornological left $A$\nbd{}module.  The natural map
  $A\hot_A W\to W$ is always injective.  The map $b'_0\colon A\hot
  W\to W$ is a bornological quotient map if and only if the map
  $A\hot_A W\to W$ induced by~$b'_0$ is a bornological isomorphism.
\end{lemma}

\begin{proof}
  Everything follows once we know that the range of $b'_1\colon A\hot
  A\hot W\to A\hot W$ is dense in the kernel of $b'_0\colon A\hot W\to
  W$.  Pick $\omega\in \Ker b'_0$.  Then there exist bounded
  complete disks $S\subseteq A$, $T\subseteq W$ such that $\omega\in
  A_S\hot W_T$.  Since $A_S$ and~$W_T$ are Banach spaces, we can find
  null sequences $(a_n)$ in~$A_S$, $(w_n)$ in~$W_T$ and $(\lambda_n)$
  in $\ell^1(\N)$ such that $\omega=\sum \lambda_n a_n\otimes w_n$
  (see~\cite{Grothendieck:Produits_Tensoriels}).  Since the set
  $\{a_n\}$ is bornologically compact in~$A$, there is a sequence
  $(u_m)$ in~$A$ such that $u_ma_n\to a_n$ for $m\to\infty$ uniformly
  for $n\in\N$.  Thus $u_m\cdot\omega\to\omega$ for $m\to\infty$.  We
  have
  $$
  b'_1(u_m\otimes\omega)=
  u_m\cdot\omega- u_m\otimes b'_0(\omega)=
  u_m\cdot\omega.
  $$
  Thus~$\omega$ is the limit of a sequence in the range of~$b'_1$.
\end{proof}

\begin{definition}  \label{def:essential_module}
  Let~$A$ be a bornological algebra with approximate identity.  A
  bornological left $A$\nbd{}module~$V$ is called \emph{essential} if
  the map $b'_0\colon A\hot V\to V$ is a bornological quotient map or,
  equivalently, $A\hot_A V\cong V$.  Essential right modules and
  bimodules are defined analogously.
\end{definition}

If~$A$ is unital, then a left $A$\nbd{}module is essential if and only
if it is unital, that is, $1_A$ acts as the identity.  The term
``essential'' is a synonym for ``non\brd{}degenerate'', which is not as
widely used for other purposes.  Grønbæk (\cite{Groenbaek:Morita})
calls such modules ``$A$\nbd{}induced''.

Let~$\ModG$ be the category of all bornological left modules over
$\CCINF(G)$.  Let $\EssG$ be its full subcategory of essential left
modules.  We write $V\in\EssG$ if~$V$ is an object of $\EssG$ and
write $f\ast v$ for $f\in\CCINF(G)$, $v\in V$, for the module
structure.

\begin{proposition}  \label{pro:essential_smooth}
  For any $V\in\EssG$ there is a natural smooth representation
  $\pi\colon G\to\Aut(V)$ such that
  $$
  f\ast v = \IN\pi(f\,dg)(v) = \int_G \pi(g,v) \cdot f(g)\,dg
  $$
  for all $f\in\CCINF(G)$, $v\in V$.  Naturality means that bounded
  module homomorphisms are $\pi$\nbd{}equivariant.
\end{proposition}

\begin{proof}
  Since~$V$ is essential, it is naturally isomorphic to the cokernel
  of the operator $b'_1\colon \CCINF(G)\hot\CCINF(G)\hot
  V\to\CCINF(G)\hot V$.  We let~$G$ act on the source and target
  of~$b'_1$ by the left regular representation on the first tensor
  factor.  This representation is smooth by
  Lemma~\ref{lem:CCINF_smooth} and~$b'_1$ is $G$\nbd{}equivariant.
  Therefore, its cokernel~$V$ carries a representation $\pi\colon
  G\to\Aut(V)$, which is smooth by Lemma~\ref{lem:smooth_hereditary}.
  It is trivial to check $\IN\pi(f_1\,dg)(f_2\ast v)=f_1\ast f_2\ast
  v$.  Since~$V$ is essential, this implies $\IN\pi(f\,dg)(v)=f\ast v$
  for all $f\in\CCINF(G)$, $v\in V$.  The construction of~$\pi$ is
  evidently natural.
\end{proof}

\subsection{Representations as modules over convolution algebras}
\label{sec:integrated_form}

We have seen how an essential module over $\CCINF(G)$ can be turned
into a smooth representation of~$G$.  Conversely, we now turn a
continuous representation $\pi\colon G\to\Aut(V)$ into a module over
$\CCINF(G)$.  Continuity implies that $Wf(g)\defeq \pi_g f(g)$ defines
a bounded linear operator from $\CCINF(G,V)$ to $L^1(G,V)\defeq
L^1(G)\hot V$, where $L^1(G)$ carries the von Neumann bornology.  We
remark without proof that the converse implication also holds: if~$W$
is a bounded linear map $\CCINF(G,V)\to L^1(G,V)$, then~$\pi$ is
already continuous.  If~$\pi$ is continuous, then
$$
\IN\pi(f\otimes v)\defeq \int_G \pi_g(v)\cdot f(g)\,dg
$$
defines a bounded linear map from $\CCINF(G,V)\cong\CCINF(G)\hot V$
to~$V$.  By adjoint associativity we obtain a bounded linear map
$\IN\pi\colon \CCINF(G)\to\End(V)$.  It is straightforward to check
that this is an algebra homomorphism, so that~$V$ becomes a module
over $\CCINF(G)$.  A morphism in $\NRepG$ between continuous
representations is a $\CCINF(G)$\nbd{}module homomorphism as well.  That
is, we have a functor from the category of continuous representations
of~$G$ to $\ModG$.

\begin{proposition}  \label{pro:section_integrated_form}
  Let $\pi\colon G\to\Aut(V)$ be a continuous representation.  Then
  the following assertions are equivalent:
  \begin{enumerate}[(i)]
  \item $\pi$ is a smooth representation, that is, the adjoint
    of~$\pi$ is a bounded linear map $V\to\CINF(G,V)$;
    
  \item the map $\IN\pi\colon \CCINF(G,V)\to V$ has a bounded
    linear right section, that is, there is a bounded linear map
    $\sigma\colon V\to\CCINF(G,V)$ such that
    $\IN\pi\circ\sigma=\ID_V$;

  \item $V$ is an essential module over $\CCINF(G)$, that is, the map
    $\IN\pi\colon \CCINF(G,V)\to V$ is a bornological quotient map.

  \end{enumerate}
  If~$\pi$ is smooth, then the section~$\sigma$ in~(ii) can be
  constructed explicitly as follows.  Choose $\phi\in\CCINF(G)$ with
  $\int_G \phi(g)\,dg=1$ and define
  $$
  \sigma_\phi\colon V\to \CCINF(G,V),
  \qquad
  \sigma_\phi(v)(g) \defeq \phi(g) \pi(g^{-1},v).
  $$
  If $H\subseteq G$ is compact, the section~$\sigma$ in~(ii) can be
  chosen $H$\nbd{}equivariant.
\end{proposition}

\begin{proof}
  If~$\pi$ is smooth, then the formula for~$\sigma_\phi$ defines a
  bounded linear map into $\CCINF(G,V)$ by
  Lemma~\ref{lem:smoothness_local}.  A trivial computation shows that
  $\sigma_\phi$ is a section for~$\IN\pi$.  Thus (i) implies~(ii).  If
  $H\subseteq G$ is compact, we can choose~$\phi$
  left-$H$\brd{}invariant.  Then the operator~$\sigma_\phi$ is
  $H$\nbd{}equivariant.  The implication (ii)$\Longrightarrow$(iii) is
  trivial.  Suppose~(iii).  The map $\IN\pi\colon \CCINF(G,V)\to V$ is
  equivariant with respect to the left regular representation of~$G$
  on $\CCINF(G,V)$.  The latter is smooth by
  Lemma~\ref{lem:CCINF_smooth}.  Thus~$\pi$ is a quotient of a smooth
  representation.  Lemma~\ref{lem:smooth_hereditary} shows that~$\pi$
  is smooth.
\end{proof}

\begin{theorem}  \label{the:essential_smooth}
  Let~$G$ be a locally compact group.  Then the categories of smooth
  representations and of essential modules are isomorphic.  The
  isomorphism sends a representation $\pi\colon G\to\Aut(V)$ to its
  integrated form $\IN\pi\colon \CCINF(G)\to\End(V)$.  In particular,
  $\pi$ is smooth if and only if~$\IN\pi$ is essential.
\end{theorem}

\begin{proof}
  The two constructions in Propositions \ref{pro:essential_smooth}
  and~\ref{pro:section_integrated_form} are clearly inverse to each
  other.  They provide the desired isomorphism of categories.
\end{proof}

\subsection{Constructions with modules and homological algebra}
\label{sec:constructions_mod}

Most functors between module categories are special cases of two
constructions: the balanced tensor product and the $\Hom$ functor.
Let~$W$ be a $B$\nbd{}$A$\brd{}bimodule.  Then we have a functor
$W\hot_A\blank$ from left $A$\nbd{}modules to left $B$\nbd{}modules and a
functor $\Hom_B(W,\blank)$ from left $B$\nbd{}modules to left
$A$\nbd{}modules.  The left $A$\nbd{}module structure on $\Hom_B(W,V)$ is
given by $a\cdot L(w)\defeq L(w\cdot a)$.  These two functors are
linked by the adjoint associativity relation
\begin{equation}  \label{eq:adjoint_associativity}
  \Hom_B(W\hot_A V,X) \cong \Hom_A(V,\Hom_B(W,X)).
\end{equation}
Of course, there are similar constructions for right modules.

Let $H\subseteq G$ be a closed subgroup.  The embedding $H\subseteq G$
induces an algebra homomorphism $\CINF'(H)\to\CINF'(G)$.  Embedding
$\CCINF(H)\subseteq \CINF'(H)$ as usual, using a left Haar
measure~$d_Hh$ on~$H$, we obtain an algebra homomorphism
$\CCINF(H)\to\CINF'(G)$.  This does not suffice to define a
restriction functor $\ModG\to\ModH$.  However, we can view $\CCINF(G)$
as a bimodule over $\CCINF(H)$ on the left and $\CCINF(G)$ on the
right by $f_0\ast f_1\ast f_2\defeq (f_0\,d_Hh)\ast f_1\ast f_2$ for
$f_0\in\CCINF(H)$, $f_1,f_2\in\CCINF(G)$.  This yields two functors
\begin{alignat*}{2}
  \smooth_G^H &\colon \ModG\to\ModH,
  &\qquad \smooth_G^H(V) &\defeq \CCINF(G)\hot_{\CCINF(G)} V,
  \\
  \nci_H^G &\colon \ModH\to\ModG,
  &\qquad \nci_H^G(V) &\defeq \Hom_{\CCINF(H)}(\CCINF(G),V),
\end{alignat*}
called \emph{(smooth) restriction functor} and \emph{(rough) induction
  functor}, respectively.  An analogous formula allows us to view
$\CCINF(G)$ as a bimodule over $\CCINF(G)$ on the left and $\CCINF(H)$
on the right.  This yields two functors
\begin{alignat*}{2}
  \ci_H^G &\colon \ModH\to\ModG,
  &\qquad \ci_H^G(V) &\defeq \CCINF(G)\hot_{\CCINF(H)} V,
  \\
  \rough_G^H &\colon \ModG\to\ModH,
  &\qquad \rough_G^H(V) &\defeq \Hom_{\CCINF(G)}(\CCINF(G),V),
\end{alignat*}
called \emph{(smooth) compact induction functor} and \emph{rough
  restriction functor}, respectively.  Finally, we define
\begin{alignat*}{2}
  \smooth &\defeq \smooth_G^G= \ci_G^G\colon \ModG\to\ModG,
  &\qquad
  \smooth(V) &\defeq \CCINF(G)\hot_{\CCINF(G)} V,
  \\
  \rough &\defeq \rough_G^G= \nci_G^G\colon \ModG\to\ModG,
  &\qquad
  \rough(V) &\defeq \Hom(\CCINF(G),V),
\end{alignat*}
the \emph{smoothening} and \emph{roughening} functors.

Our treatment of the compact induction functor as a tensor product is
analogous to Marc Rieffel's approach to induced representations
(\cite{Rieffel:Induced}).  The Banach algebra variant of Rieffel's
theory by Niels Grønbæk is even closer to our setup
(\cites{Groenbaek:Morita, Groenbaek:Imprimitivity}).  The only
difference is that Grønbæk works with $L^1(G)$ instead of $\CCINF(G)$.

The following theorem shows that the smoothening deserves its name.
We use the natural map $\smooth(V)\to V$ induced by $b'_0(f\otimes
v)\defeq f\ast v$.

\begin{theorem}  \label{the:smoothening}
  The natural map $\smooth(V)\to V$ is always injective and an
  isomorphism if and only if $V\in\EssG$.  The smoothening is an
  idempotent functor on~$\ModG$ whose range is~$\EssG$.  As a functor
  $\ModG\to\EssG$ it is left adjoint to the embedding
  $\EssG\to\ModG$.  Let $\pi\colon G\to\Aut(V)$ be a continuous
  representation of~$G$.  Then the smoothenings of~$G$ as a module and
  as a representation agree.
\end{theorem}

\begin{proof}
  We know from Lemma~\ref{lem:approxid_hota} that the map
  $\smooth(V)\to V$ is always injective and an isomorphism if and only
  if~$V$ is essential.  Since the left regular representation on
  $\CCINF(G)$ is smooth, $\CCINF(G)$ is an essential left module over
  itself by Theorem~\ref{the:essential_smooth}.  That is,
  $\CCINF(G)\hot_{\CCINF(G)}\CCINF(G)\cong\CCINF(G)$.  Since the
  balanced tensor product is associative, we obtain
  $\smooth^2=\smooth$.  Since $S(V)\cong V$ if and only if
  $V\in\EssG$, the range of~$\smooth$ is~$\EssG$.

  Let~$W$ be an essential module.  Since the map $\smooth(V)\to V$ is
  always injective, the induced map $\Hom(W,\smooth(V))\to\Hom(W,V)$
  is injective.  Any bounded module homomorphism $W\to V$ restricts to
  a bounded module homomorphism $W=\smooth(W)\to\smooth(V)$, so that
  the map $\Hom(W,\smooth(V))\to\Hom(W,V)$ is also surjective.  This
  means that the embedding and smoothening functors are adjoint.
  
  Let~$\pi$ be a continuous representation.  Let $V_0$ and~$V_1$ be
  the smoothenings of~$V$ as a representation and as a module,
  respectively.  The natural maps $V_0\to V$ and $V_1\to V$ are both
  injective.  Since~$V_1$ is an essential module, it is a smooth
  representation of~$G$ as well.  Hence the map $V_0\to V$ factors
  through $V_0\to V_1$ by the universal property of the smoothening.
  Similarly, since~$V_0$ is an essential module, the map $V_1\to V$
  factors through $V_1\to V_0$.  Both maps $V_0\to V_1$ and $V_1\to
  V_0$ are injective and bounded, hence bornological isomorphisms.
\end{proof}

Equation~\eqref{eq:adjoint_associativity} specializes to natural
isomorphisms
\begin{align}
  \label{eq:adjoint_ci_rough}
  \Hom_{\CCINF(G)}(\ci_H^G(V),W) &\cong
  \Hom_{\CCINF(H)}(V,\rough_G^H(W)),
  \\
  \label{eq:adjoint_smooth_i}
  \Hom_{\CCINF(H)}(\smooth_G^H(V),W) &\cong
  \Hom_{\CCINF(G)}(V,\nci_H^G(W)).
\end{align}
That is, compact induction is left adjoint to rough restriction and
rough induction is right adjoint to smooth restriction.

Especially, $\smooth$ is left adjoint to $\rough$.  Being adjoint to
an idempotent functor, $\rough$ is idempotent as well.  Thus~$\rough$
is a projection onto a subcategory of~$\ModG$.  We may call these
modules \emph{rough}.  They are usually not smooth, but if~$G$ is a
Lie group they are differentiable by
Theorem~\ref{the:explicit_criterion_differentiable} because they are
evidently modules over $\CINF'(G)$.  We have
$\rough\circ\smooth\cong\rough$ because
\begin{multline*}
  \Hom_{\CCINF(G)}(V,\rough\circ\smooth(W))
  \cong
  \Hom_{\CCINF(G)}(\smooth(V),\smooth(W))
  \\ \cong
  \Hom_{\CCINF(G)}(\smooth(V),W)
  \cong
  \Hom_{\CCINF(G)}(V,\rough(W))
\end{multline*}
for all $V,W\in\ModG$.  We will prove shortly that
$\smooth\circ\rough\cong\smooth$.  Summarizing, we have
\begin{equation}  \label{eq:smooth_rough}
  \smooth\circ\smooth \cong\smooth,
  \quad
  \smooth\circ\rough  \cong\smooth,
  \quad
  \rough\circ\smooth  \cong\rough,
  \quad
  \rough\circ\rough   \cong\rough.
\end{equation}
The natural map $V\to\rough(V)$ is injective if and only if no
non-zero vector $v\in V$ satisfies $f\ast v=0$ for all
$f\in\CCINF(G)$.  Let us restrict attention to this class of modules.
Then the natural maps $\smooth(V)\to V\to\rough(V)$ are injective.  If
we have injective maps $\smooth(V)\to W\to\rough(V)$, then
$\smooth(V)=\smooth(W)$ because already $\smooth\rough(V)=\smooth(V)$
and the smoothening preserves monomorphisms.  Conversely, if
$\smooth(W)\cong\smooth(V)$, then
$\rough(W)\cong\rough\smooth(W)\cong\rough\smooth(V)\cong\rough(V)$ as
well, so that we have injective maps $\smooth(V)\to W\to\rough(V)$.
This means that a module~$W$ satisfies $\smooth(W)=\smooth(V)$ if and
only if it lies between $\smooth(V)$ and $\rough(V)$.

In the following we tacitly identify~$\EssG$ with~$\RepG$ using
Theorem~\ref{the:essential_smooth}.  If we have to view a smooth
representation as a right module, we always use the antipode
$\tilde{f}^{(1)}$ defined in~\eqref{eq:CCINF_antipode} to turn a left
into a right module.

Since $\smooth(V)=V$ for $V\in\EssG$, we have
$\smooth_G^H|_{\EssG}\cong\Res_G^H$.  The universal property of the
smoothening and~\eqref{eq:adjoint_smooth_i} imply that
$\smooth\circ\nci_H^G(W)\colon \EssH\to\EssG$ is right adjoint to
$\Res_G^H$.  This means that
\begin{equation}
  \label{eq:nci_Ind}
  \smooth\circ\nci_H^G \cong \Ind_H^G.
\end{equation}
Since $\Ind_G^G$ is the identical functor, we get the relation
$\smooth\circ\rough=\smooth$ claimed in~\eqref{eq:smooth_rough}.  The
relationship between $\ci_H^G$ and $\cInd_H^G$ is more complicated.
Before we discuss it we need some other useful results.

Let $X$ and~$Y$ be a right and left module over $\CCINF(G)$ and
let~$W$ be a bornological vector space.  Then $\Hom(X,W)$ is a left
module over $\CCINF(G)$ in a canonical way
and~\eqref{eq:adjoint_associativity} yields
\begin{equation}
  \label{eq:G_hot_adjoint}
  \Hom(X\hot_{\CCINF(G)} Y,W) \cong \Hom_{\CCINF(G)}(Y,\Hom(X,W)).
\end{equation}
Let $\C(1)$ be the trivial representation of~$G$ on~$\C$ viewed as a
right module over $\CCINF(G)$.  The space $\C(1)\hot_{\CCINF(G)} Y$ is
called the \emph{coinvariant space} of~$Y$.  If~$Y$ is a smooth
representation viewed as a left module over $\CCINF(G)$ and $W=\C$,
then~\eqref{eq:G_hot_adjoint} asserts that the dual space of the
coinvariant space of~$Y$ is the space of $G$\nbd{}invariant linear
functionals on~$Y$.

Let $X,Y,Z$ be smooth representations of~$G$.  We let~$G$ act on
$\Hom(Y,Z)$ by the conjugation action $(g\cdot l)(y)\defeq
g\cdot l(g^{-1}y)$ and on $X\hot Y$ by the diagonal action $g\cdot
(x\otimes y)\defeq gx\otimes gy$.  These two constructions are adjoint
in the sense that
\begin{equation}
  \label{eq:adjoint_diagonal_conjugation}
  \Hom_G(X,\smooth \Hom(Y,Z)) \cong
  \Hom_G(X,\Hom(Y,Z)) \cong \Hom_G(X\hot Y,Z).
\end{equation}
The first isomorphism is the universal property of the smoothening.
The second is proved by identifying both sides with the space of
bilinear maps $l\colon X\times Y\to Z$ that satisfy the equivariance
condition $l(gx,gy)=gl(x,y)$.  If we let $X\defeq\C(1)$ be the trivial
representation of~$G$ on~$\C$, we have $\C(1)\hot Y\cong Y$ and
\begin{equation}
  \label{eq:G_hom_invariants}
  \Hom_G(\C(1),\smooth \Hom(Y,Z)) \cong
  \Hom_G(Y,Z).
\end{equation}
Next we claim that
\begin{equation}
  \label{eq:G_hot_coinvariants}
  \C(1)\hot_{\CCINF(G)} (Y\hot Z) \cong
  Y\hot_{\CCINF(G)} Z,
\end{equation}
where we view $\C(1)$ and~$Y$ as right modules over $\CCINF(G)$.
Equation~\eqref{eq:G_hot_coinvariants} can easily be verified
directly.  For the fun of it we use adjointness relations to
prove the equivalent assertion that $\Hom(\C(1)\hot_{\CCINF(G)} (Y\hot
Z),W)\cong\Hom(Y\hot_{\CCINF(G)} Z,W)$ for all bornological vector
spaces~$W$.  Equation~\eqref{eq:G_hot_adjoint} implies
\begin{align*}
  \Hom(Y\hot_{\CCINF(G)} Z,W) &\cong
  \Hom_G(Z,\Hom(Y,W)),
  \\
  \Hom(\C(1)\hot_{\CCINF(G)} (Y\hot Z),W) &\cong
  \Hom(Y\hot Z,\Hom(\C(1),W))\cong \Hom(Y\hot Z,W),
\end{align*}
where~$G$ acts on $\Hom(Y,W)$ by $g\cdot l(y)\defeq l(g^{-1}y)$ and
trivially on~$W$.  Since the action on $\Hom(Y,W)$ is the conjugation
action for the trivial representation on~$W$, both spaces are
isomorphic by~\eqref{eq:adjoint_diagonal_conjugation}.  This finishes
the proof of~\eqref{eq:G_hot_coinvariants}.

Now we are ready to relate the functors $\ci_H^G$ and $\cInd_H^G$.
Recall that $\modular_{G:H}$ denotes the quasi\brd{}character
$\modular_G/\modular_H\colon H\to\R\inv_+$.  For a representation
$\pi\colon H\to\Aut(V)$ we write $\modular_{G:H}\cdot V$ for the
representation $\modular_{G:H}\cdot\pi$ on~$V$.

\begin{theorem}  \label{the:ci_cInd}
  There is a natural isomorphism $\ci_H^G(V) \cong
  \cInd_H^G(\modular_{G:H}\cdot V)$ for all $V\in\EssH$.
\end{theorem}

\begin{proof}
  First we explain the source of the relative modular function in
  $\ci_H^G(V)$.  The right $\CCINF(G)$\brd{}module structure on
  $\CCINF(G)$ is the integrated form of the twisted right regular
  representation $\rho\cdot\modular_G$ because $f(g)\,d_Gg\ast
  \delta_{x^{-1}}=f(gx)\modular_G(x)\,d_Gg$.  We equip $\CCINF(G)$ and
  $\CCINF(H)$ with the canonical $\CCINF(H)$\brd{}bimodule structure.
  The restriction map $\CCINF(G)\to\CCINF(H)$ is a left module
  homomorphism, but we pick up a factor $\modular_{G:H}$ for the right
  module structure.  Therefore, it induces an $H$\nbd{}equivariant map
  $\ci_H^G(V)\to \modular_{G:H}\cdot V$ and hence a $G$\nbd{}equivariant
  map into $\Ind_H^G(\modular_{G:H}\cdot V)$.  This is the desired
  isomorphism onto $\cInd_H^G(\modular_{G:H}\cdot V)$.  We now
  construct it more explicitly.  Define
  $$
  \Phi \colon \CCINF(G,V)\to\CINF(G,V),
  \qquad
  \Phi f(g)\defeq
  \int_H h\cdot f(g^{-1}\cdot h)\modular_{G:H}(h)\,d_Hh,
  $$
  where~$d_Hh$ is a left invariant Haar measure on~$H$.  Clearly,
  $\supp\Phi f \subseteq H\cdot (\supp f)^{-1}$ is uniformly compact
  in $H\backslash G$ for $f$ in a bounded subset of $\CCINF(G,V)$.
  Moreover, $\Phi f(hg)=\modular_{G:H}(h) h\cdot \Phi f(g)$ for all
  $h\in H$, $g\in G$.  This means that the range of~$\Phi$ is
  contained in $\cInd_H^G(\modular_{G:H}V)$.  Moreover, one computes
  easily that $\Phi(f_h)=\Phi(f)$ if $f_h(g)\defeq \modular_G(h)
  h\cdot f(gh)$ for $h\in H$.  This means that~$\Phi$ is
  $H$\nbd{}invariant for the diagonal action of~$H$ on $\CCINF(G)\hot V$
  that occurs in~\eqref{eq:G_hot_coinvariants}.  Therefore, $\Phi$
  descends to a bounded linear map on $\CCINF(G)\hot_{\CCINF(H)}
  V=\ci_H^G(V)$.  Finally, $\Phi$ is $G$\nbd{}equivariant, that is,
  $\Phi\lambda_g=\rho_g\Phi$.  Summing up, we have constructed a
  natural transformation
  $$
  \Phi\colon \ci_H^G(V)\to\cInd_H^G(\modular_{G:H}\cdot V).
  $$
  
  It remains to verify that~$\Phi$ is an isomorphism for all~$V$.
  This is easy for the left regular representations on $\CCINF(H,V)$,
  where we can compute both sides explicitly.  Any essential module
  over $\CCINF(H)$ is the cokernel of a map $b'_1\colon \CCINF(H\times
  H,V)\to\CCINF(H,V)$ between left regular modules.  The functor
  $\ci_H^G$ preserves cokernels because it has a right adjoint.  The
  functor $\cInd_H^G$ also preserves cokernels by
  Proposition~\ref{pro:cInd_concrete}.  Hence~$\Phi$ is an isomorphism
  for all~$V$.
\end{proof}

\begin{corollary}  \label{cor:induction_cocompact}
  If $H\subseteq G$ is cocompact, then there is a natural isomorphism
  $$
  \smooth\circ\nci_H^G(\modular_{G:H}\cdot V)\cong\ci_H^G(V).
  $$
\end{corollary}

\begin{proof}
  It is clear from the definition that $\cInd_H^G=\Ind_H^G$ in this
  case.  Hence the assertion follows from Theorem~\ref{the:ci_cInd}
  and~\eqref{eq:nci_Ind}.
\end{proof}

We continue with some further properties of our functors.  Let
$L\subseteq H\subseteq G$.  Since the right $\CCINF(H)$\brd{}module
structure on $\CCINF(G)$ comes from a smooth representation, we have
$\CCINF(G)\hot_{\CCINF(H)}\CCINF(H)\cong \CCINF(G)$ and hence
$$
\ci_H^G\circ\ci_L^H(V) =
\CCINF(G)\hot_{\CCINF(H)}\CCINF(H)\hot_{\CCINF(L)} V \cong
\CCINF(G)\hot_{\CCINF(L)} V =
\ci_L^G(V).
$$
The assertion $\smooth_H^L\circ\smooth_G^H=\smooth_G^L$ is proved
similarly.  By adjointness we also obtain
$\nci_H^G\circ\nci_L^H=\nci_L^G$ and
$\rough_H^L\circ\rough_G^H=\rough_G^L$.  We evidently have
$\Res_H^L\Res_G^H=\Res_G^L$ and hence $\Ind_H^G\circ\Ind_L^H=\Ind_L^G$
by adjointness.  As special cases we note that
\begin{equation}
  \label{eq:induce_smooth_rough}
  \rough\circ \nci_H^G=\nci_H^G=\nci_H^G\circ\rough,
  \qquad
  \smooth\circ\ci_H^G=\ci_H^G\circ\smooth=\ci_H^G.
\end{equation}
Together with~\eqref{eq:smooth_rough}, we obtain further relations
like $\nci_H^G\circ\smooth=\nci_H^G$ and $\ci_H^G\circ\rough=\ci_H^G$.

Let $V$ and~$W$ be a right and a left module over $\CCINF(G)$ and
$\CCINF(H)$, respectively.  Then we trivially have
\begin{equation}
  \label{eq:hot_ci_smooth}
  V\hot_{\CCINF(G)} \ci_H^G(W) \cong
  V\hot_{\CCINF(G)} \CCINF(G) \hot_{\CCINF(H)} W \cong
  \smooth_G^H V\hot_{\CCINF(H)} W.
\end{equation}
Let~$X$ be a bornological vector space, equip $\Hom(W,X)$ with the
canonical right module structure.  Then we have canonical isomorphisms
\begin{multline*}
  \Hom_{\CCINF(G)}(V,\nci_H^G\Hom(W,X))
  \cong
  \Hom_{\CCINF(H)}(\smooth_G^H V,\Hom(W,X))
  \\ \cong
  \Hom(\smooth_G^H V \hot_{\CCINF(H)} W,X)
  \cong
  \Hom(V\hot_{\CCINF(G)} \ci_H^G W,X)
  \\ \cong
  \Hom_{\CCINF(G)}(V,\Hom(\ci_H^G W,X)).
\end{multline*}
Since~$V$ is arbitrary, we conclude that
\begin{equation}
  \label{eq:induce_dual}
  \nci_H^G\Hom(W,X)\cong \Hom(\ci_H^G W,X).
\end{equation}
as left modules over $\CCINF(G)$.  Here~$W$ is a right module over
$\CCINF(H)$ and~$X$ is a bornological vector space.  For $X=\C$ this
is an assertion about induction of dual spaces.  The smoothening of
the dual is the \emph{contragradient} representation~$\tilde{W}$.
Equation~\eqref{eq:induce_dual} implies
\begin{equation}
  \label{eq:induce_contragradient}
  \Ind_H^G \tilde{W}\cong
  \bigl(\cInd_H^G(\modular_{G:H}\cdot W)\bigr)\sptilde.
\end{equation}
The analogous statements
\begin{equation}
  \label{eq:restrict_dual}
  \rough_G^H \Hom(W,X) \cong
  \Hom(\smooth_G^H W,X),
  \qquad
  \rough \Hom(W,X)\cong\Hom(\smooth W,X),
\end{equation}
about restriction follow easily
from~\eqref{eq:adjoint_associativity}.

Finally, we do some homological algebra and begin by recalling a few
standard notions.  Let~$A^+$ be the augmented unital algebra obtained
by adjoining a unit element to a bornological algebra~$A$.  The
category of left modules over~$A$ is isomorphic to the category of
unital left modules over~$A^+$.  Hence the correct definition of a
\emph{free left module} over~$A$ is $A^+\hot V$ with the evident left
module structure over~$A$.  Similar remarks apply to right modules and
bimodules.  The free module has the universal property that bounded
module homomorphisms $A^+\hot V\to W$ correspond bijectively to
bounded linear maps $V\to W$.  As a consequence, free modules are
projective for linearly split extensions.  In the following we say
that a module is \emph{relatively projective} if it is projective for
this class of extensions.  In general, the modules $A\hot V$ need not
be relatively projective.

\begin{proposition}  \label{pro:CCINF_projective}
  Let $H\subseteq G$.  Then $\CCINF(G)$ is relatively projective as a
  left or right module over $\CCINF(H)$.
\end{proposition}

\begin{proof}
  It suffices to prove that $\CCINF(G)$ is projective as a left module
  over $\CCINF(H)$.  We are going to construct a bounded
  $\CCINF(H)$\brd{}linear section~$\sigma$ for the convolution map
  $$
  \mu\colon \CCINF(H\times G)\cong
  \CCINF(H)\hot\CCINF(G)\to\CCINF(G),
  \qquad
  \mu f(g)\defeq \int_H f(h,h^{-1}g)\,dh.
  $$
  Let~$\mu^+$ be the extension of~$\mu$ to $\CCINF(H)^+\hot\CCINF(G)$,
  then $\mu^+\circ\sigma=\ID$ as well.  Thus $\CCINF(G)$ is relatively
  projective as a retract of the free module
  $\CCINF(H)^+\hot\CCINF(G)$.  The map $\sigma$ is defined by $\sigma
  f(h,g)\defeq f(hg)\cdot \phi(g)$ for some function
  $\phi\in\CINF(G)$.  This defines a map to $\CCINF(H\times G)$ if
  $\supp\phi \cap H\cdot L$ is compact for all compact $L\subseteq G$.
  It is a section for~$\mu$ if and only if $\int_H
  \phi(h^{-1}g)\,d_Hh=1$ for all $g\in G$.  Functions~$\phi$ with
  these properties clearly exist.
\end{proof}

\begin{theorem}  \label{the:smoothen_exact}
  Let $H\subseteq G$ be a closed subgroup.   The functors
  $\ci_H^G$ and $\smooth_H^G$ preserve bornological extensions,
  locally linearly split extensions, linearly split extensions and
  injectivity of morphisms.  They commute with arbitrary direct
  limits.  They map relatively projective objects to relatively
  projective objects.  In particular, all this applies to the
  smoothening functor.
  
  The functors $\nci_H^G$ and $\rough_H^G$ preserve linearly split
  extensions and injectivity of morphisms.  They commute with
  arbitrary inverse limits.  They map relatively injective objects to
  relatively injective objects.  In particular, all this applies to
  the roughening functor.
\end{theorem}

\begin{proof}
  For the exactness assertions we can forget the module structure on
  $\ci_H^G(V)$ and $\smooth_G^H(V)$ and view these spaces just as
  bornological vector spaces.  Thus the exactness assertions about
  $\smooth_G^H$ follow from the corresponding statements about
  $\ci_G^G$.  Proposition~\ref{pro:CCINF_projective} implies that the
  functor $\ci_H^G$ is a retract of the functor $V\mapsto
  \CCINF(G)\hot V\cong \CCINF(G,V)$.  Hence it inherits the properties
  of the latter functor listed in Proposition~\ref{pro:CCINF_exact}.
  Since $\ci_H^G$ and $\smooth_G^H$ have right adjoints, they commute
  with direct limits.  Furthermore, the assertion that $\ci_H^G$
  preserves relative projectivity is equivalent to the statement that
  its right adjoint functor $\rough_H^G$ is exact for linearly split
  extensions.  This follows from
  Proposition~\ref{pro:CCINF_projective}.  It is evident that
  $\nci_H^G$ and $\rough_H^G$ preserve injectivity of morphisms.
  Since they have left adjoint functors, they commute with inverse
  limits.  Since their left adjoints are exact for linearly split
  extensions, they preserve relatively injective objects.
\end{proof}

\begin{theorem}  \label{the:essential_extension}
  Let $K\into E\prto Q$ be a bornological extension in~$\ModG$.  Then
  $E\in\EssG$ if and only if both $K\in\EssG$ and~$Q\in\EssG$.
\end{theorem}

\begin{proof}
  Let $K',E',Q'$ be the smoothenings of $K,E,Q$.  Consider the diagram
  $$
  \xymatrix{
    {K'\;} \ar[d] \ar@{>->}[r] & {E'} \ar[d] \ar@{->>}[r] & {Q'} \ar[d] \\
    {K\;} \ar@{>->}[r] & {E} \ar@{->>}[r] & {Q.} \\
  }
  $$
  Both rows are bornological extensions by
  Theorem~\ref{the:smoothen_exact}.  If $K$ and~$Q$ are essential,
  then the vertical arrows $K'\to K$ and $Q'\to Q$ are bornological
  isomorphisms.  This implies that the middle arrow is a bornological
  isomorphism by the Five Lemma.  The validity of the Five Lemma for
  bornological vector spaces can be proved directly.  It also follows
  easily from the observation that the category of bornological vector
  spaces with the class of bornological extensions is an exact
  category in the sense of Daniel Quillen (see~\cites{Quillen:Higher,
  Prosmans-Schneiders:Bornological}).  Hence~$E$ is essential if both
  $K$ and~$Q$ are essential.  Conversely, if~$E$ is essential, then
  the module action $\CCINF(G)\hot Q\to Q$ is a bornological quotient
  map, so that~$Q$ is essential.  Another application of the Five
  Lemma shows that~$K$ is essential as well.
\end{proof}

We have seen in Section~\ref{sec:smooth_permanence} that the class of
smooth representations of~$G$ is hereditary for subrepresentations and
quotient representations, but not for extensions in general.  We have
to assume the representation on~$E$ to be continuous.  Then we can use
Theorem~\ref{the:essential_extension} to obtain the smoothness
of~$E$.

\begin{theorem}  \label{the:enough_pro_inj}
  The category $\EssG\cong\RepG$ has enough relatively projective and
  injective objects.
  
  The functor $\Ind_H^G\colon \RepH\to\RepG$ is exact for linearly
  split extensions.  It preserves monomorphisms and relatively
  injective objects.  It commutes with inverse limits in these
  subcategories (they differ from those in the larger categories
  $\NRepG$ or $\ModG$!).
  
  The functors $\cInd_H^G$ and $\Res_G^H$ are exact for any class of
  extensions and preserve monomorphisms and relatively projective
  objects.  They commute with direct limits.
\end{theorem}

\begin{proof}
  The exactness assertions about $\Res_G^H$ are trivial.  The
  exactness properties of $\Ind_H^G\cong\smooth\circ\nci_H^G$ follow
  immediately from those of $\smooth$ and $\nci_H^G$.  Since
  $\Res_H^G$ and $\Ind_H^G$ are adjoint, the first preserves direct
  and the latter preserves inverse limits.  The exactness properties
  imply that $\Ind_H^G$ and $\Res_G^H$ preserve relatively injective
  and projective objects, respectively.  The assertions about
  $\cInd_H^G$ follow immediately from the corresponding properties of
  $\ci_H^G$ and Theorem~\ref{the:ci_cInd}.  For the trivial group~$E$,
  linearly split extensions are already direct sum extensions.  Thus
  any object is relatively injective and projective.  By
  Theorem~\ref{the:smoothen_exact} we obtain that
  $\ci_E^G(V)=\CCINF(G,V)$ is relatively projective and
  $\Ind_E^G(V)=\smooth \CINF(G,V)$ is relatively injective.  If~$V$ is
  an arbitrary smooth representation, then we have a linearly split
  surjection $\CCINF(G,V)\to V$ by
  Proposition~\ref{pro:section_integrated_form} and a linearly split
  injection $V\to \smooth \CINF(G,V)$.
\end{proof}

Thus we can derive functors on the category of smooth representations
using relatively projective and injective resolutions.  Let us write
$\Left_\ast F$ and $\Right^\ast F$, $\ast\in\N$, for the left and
right derived functors of a functor~$F$ from $\RepG$ to some additive
category.  The left derived functors of $V\hot_{\CCINF(G)}\blank$ are
denoted $\Tor_\ast^G(V,W)$, the right derived functors of
$\Hom_G(V,\blank)$ are denoted $\Ext^\ast_G(V,W)$.  If we take~$V$ to
be the trivial representation on~$\C$, we obtain group homology and
cohomology, denoted $\Ho_\ast(G,V)$ and $\Ho^\ast(G,V)$, respectively.

The general machinery of derived functors yields the following
results.  Since the compact induction functor is exact and preserves
relatively projective objects, we have
$\Left_\ast(F\circ\ci_H^G)=(\Left_\ast F)\circ \ci_H^G$.  Since the
induction functor $\Ind_H^G$ is exact and preserves relatively
injective objects, we have $\Right^\ast(F\circ\Ind_H^G)=(\Right^\ast
F)\circ\Ind_H^G$.  Therefore, the adjointness of restriction and
induction and~\eqref{eq:hot_ci_smooth} imply
\begin{align}
  \label{eq:induce_ext}
  \Ext^\ast_G(V,\Ind_H^G(W)) &\cong
  \Ext^\ast_H(\Res_G^H V,W),
  \\
  \label{eq:induce_tor}
  \Tor_\ast^G(V,\cInd_H^G(\modular_{G:H}\cdot W)) &\cong
  \Tor_\ast^H(\Res_G^H V,W),
  \\
  \label{eq:Shapiro_cohomology}
  \Ho^\ast(G,\Ind_H^G(W))
  &\cong \Ho^\ast(H,W),
  \\
  \label{eq:Shapiro_homology}
  \Ho_\ast(G,\cInd_H^G(\modular_{G:H}\cdot W))
  &\cong \Ho_\ast(H,W).
\end{align}

The functors $W\mapsto V\hot W$ with diagonal action and $W\mapsto
\Hom(V,W)$ with conjugation action are evidently exact for linearly
split extensions.  Since they are adjoint
by~\eqref{eq:adjoint_diagonal_conjugation}, the first preserves
relative projectivity and the second preserves relative injectivity.
Reasoning as above \eqref{eq:G_hom_invariants}
and~\eqref{eq:G_hot_coinvariants} imply
\begin{align}
  \label{eq:ext_univariant}
  \Ext^\ast_G(V,W) & \cong
  \Ext^\ast_G(\C(1),\smooth \Hom(V,W)) =
  \Ho^\ast(G,\smooth \Hom(V,W)),
  \\
  \label{eq:tor_univariant}
  \Tor_\ast^G(V,W) & \cong
  \Tor_\ast^G(\C(1),V\hot W) =
  \Ho_\ast(G, V\hot W).
\end{align}
That is, group homology and cohomology already determine the bivariant
homology and cohomology theories.

\subsection{The Gårding subspace}
\label{sec:Garding_subspace}

The smoothening for modules is closely related to the Gårding
subspace.  Let~$V$ be a continuous representation of a locally compact
group on a bornological vector space.  The \emph{Gårding subspace}
of~$V$ is defined as the linear subspace spanned by $\IN\pi(f)(v)$
with $f\in\CCINF(G)$, $v\in V$.  This is the image of the uncompleted
tensor product $\CCINF(G)\otimes V$ in~$V$.  In contrast, $\smooth(V)$
is the image of the completed tensor product $\CCINF(G)\hot V$.  It
seems that everything that can be done with the Gårding subspace can
also be done with $\CCINF(G)\hot_{\CCINF(G)} V$.  However, it is
actually true that the Gårding subspace is always equal to
$\smooth(V)$.  This is proved by Jacques Dixmier and Paul Malliavin
in~\cite{Dixmier-Malliavin:Factorisations} for Lie group
representations on Fréchet spaces.  The same argument actually works
in much greater generality:

\begin{theorem}  \label{the:Garding_subspace}
  Let $\pi\colon G\to\Aut(V)$ be a continuous representation of a
  locally compact group~$G$ on a bornological vector space~$V$.  The
  Gårding subspace of~$V$ is equal to $\smooth(V)$.  Especially, any
  element of $\CCINF(G)$ is a finite linear combination of products
  $f_1\ast f_2$ with $f_1,f_2\in\CCINF(G)$.
\end{theorem}

\begin{proof}
  We may assume that the representation~$V$ is already smooth because
  we only make the problem more difficult if we shrink~$V$ to
  $\smooth(V)$.  Any $v\in V$ already belongs to $V^k$ for some
  smooth compact subgroup $k\subseteq G$.  We can replace the
  representation of~$G$ on~$V$ by the smooth representation of the Lie
  group $N_G(k)/k$ on~$V^k$.  Thus we may assume~$G$ to be a Lie group
  without loss of generality.  The class of smooth representations for
  which the theorem holds is evidently closed under inductive limits
  and under quotients.  If~$V$ is a smooth representation, then it is
  a quotient of the left regular representation on $\CCINF(G,V)$.  The
  latter is the inductive limit of the left regular representations on
  $\CCINF(G,V_T)$ for the small complete disks $T\subseteq V$.
  Hence it suffices to prove the assertion for the left regular
  representation on $\CCINF(G,V_T)$ for a Banach space~$V_T$.  This
  case can be dealt with by literally the same argument that Jacques
  Dixmier and Paul Malliavin use
  in~\cite{Dixmier-Malliavin:Factorisations} to prove that the Gårding
  subspace of $\CCINF(G)$ is $\CCINF(G)$.
\end{proof}

\section{The center of the category of smooth representations}
\label{sec:CCINF_Mult}

\begin{definition}  \label{def:multipliers}
  Let~$A$ be a bornological algebra with the property that $A\cdot A$
  spans a dense subspace of~$A$.
  
  Let $\LMult(A)$ and $\RMult(A)^\op$ be the algebras of bounded right
  and left module homomorphisms $A\to A$, equipped with the
  equibounded bornology.  These are the \emph{left and right
  multiplier algebras} of~$A$.  By convention, the multiplication in
  $\RMult(A)$ is the opposite of the composition of operators.  The
  \emph{(two-sided) multiplier algebra} $\Mult(A)$ of~$A$ is the
  algebra of pairs $(l,r)$ of a left and a right multiplier such that
  $a\cdot (l\cdot b)=(a\cdot r)\cdot b$ for all $a,b\in A$.
\end{definition}

All three multiplier algebras are unital bornological algebras and
there are obvious bounded algebra homomorphisms from~$A$ into them.
We claim that~$A$ is a bornological unital
$\LMult(A)$\brd{}$\RMult(A)$\brd{}bimodule.  The only point that is
not obvious is that $(l\cdot a)\cdot r=l\cdot(a\cdot r)$ for all $a\in
A$, $l\in\LMult(A)$, $r\in\RMult(A)$.  If $a=bc$ with $b,c\in A$, then
$(l\cdot bc)\cdot r=(lb)\cdot (cr)=l\cdot(bc\cdot r)$.  The claim
follows because the linear span of elements of the form~$bc$ is dense
in~$A$.

We denote the center of an algebra~$A$ by $\Center(A)$.  A left
multiplier~$l$ of~$A$ is called \emph{central} if $a\cdot l\cdot
b=l\cdot a\cdot b$ for all $a,b\in A$.  That is, the pair $(l,l)$ is a
two-sided multiplier of~$A$.  Since we know that left and right
multipliers commute with each other, it follows that~$l$ commutes with
any left or right multiplier on~$A$.  Thus~$l$ belongs to the centers
of all three multiplier algebras.  Conversely, if $l$ is central, say,
in $\LMult(A)$, then it is a central multiplier in the above sense
because $A\subseteq\LMult(A)$.  As a result, the multiplier algebras
all have the same center, which consists exactly of the central
multipliers.

\begin{definition}  \label{def:center_category}
  The \emph{center} $\Center(\Cat)$ of an additive category~$\Cat$ is
  the ring of natural transformations from the identity functor
  $\ID\colon \Cat\to\Cat$ to itself.
\end{definition}

Equivalently, an element of~$\Center(\Cat)$ is a family of morphisms
$\gamma_X\colon X\to X$ for each object~$X$ of~$\Cat$ such that
$f\circ\gamma_X=\gamma_Y\circ f$ for any morphism $f\colon X\to Y$
in~$\Cat$.  The center of the category of smooth representations of a
totally disconnected group on vector spaces is studied by Joseph
Bernstein in~\cite{Bernstein:Center} and plays a crucial role in the
representation theory of reductive groups over non\brd{}Archimedean
local fields.

\begin{lemma}  \label{lem:center_approxid}
  Let~$A$ be a bornological algebra with an approximate identity.
  Suppose that $A\hot_A A\cong A$.  Then the center of the category of
  essential $A$\nbd{}modules is naturally isomorphic to the algebra of
  central multipliers of~$A$.
\end{lemma}

\begin{proof}
  Let~$\Cat$ be the category of essential bornological left
  $A$\nbd{}modules.  The center of~$\Cat$ maps into the center of the
  endomorphism ring of~$A$ because $A\in\Cat$.  By definition, this
  endomorphism ring is $\RMult(A)^\op$.  Hence its center is the
  algebra of central multipliers.  Thus we obtain a homomorphism
  $\alpha\colon \Center(\Cat)\to\Center\Mult(A)$.  We have to check
  that this map is bijective.
  
  For injectivity suppose that $\Phi\in\Center(\Cat)$ vanishes
  on~$A$.  Let $V\in\Cat$ and $v\in V$.  Then the map $a\mapsto av$ is
  a morphism $A\to V$ in~$\Cat$.  Hence $\Phi_V(av)=\Phi_A(a)v=0$.
  Since elements of the form $av$ generate~$V$, we get $\Phi_V=0$.
  Thus~$\alpha$ is injective.  For surjectivity let~$l$ be a central
  multiplier.  Since~$A$ is a bimodule over $\LMult(A)$ and~$A$, there
  is a canonical $\LMult(A)$\brd{}module structure on $A\hot_A V$, that
  is, on any essential module.  Thus~$l$ acts in a canonical way on
  any $V\in\Cat$.  Centrality implies that~$l$ acts by left module
  homomorphisms.  Thus we obtain an element of $\Center(\Cat)$.
\end{proof}

The center of the category of all modules over~$A$ is equal to the
center of~$A^+$ because modules over~$A$ are the same as essential
modules over~$A^+$.  Hence we may get a much smaller center than for
essential modules.

\begin{theorem}  \label{the:center_smoothrep}
  Let~$G$ be a locally compact group.  Then the center of the category
  of smooth representations of~$G$ is naturally isomorphic to $\Center
  \Mult(\CCINF(G))$, the algebra of central multipliers of
  $\CCINF(G)$.
\end{theorem}

\begin{proof}
  Theorem~\ref{the:essential_smooth} asserts that $\RepG$ is
  isomorphic to~$\EssG$ and hence has an isomorphic center.  We know
  that $\CCINF(G)$ satisfies the hypotheses of
  Lemma~\ref{lem:center_approxid}.  Hence $\Center(\RepG)\cong\Center
  \Mult(\CCINF(G))$.
\end{proof}

\begin{lemma}  \label{lem:CCINF_mult_distribution}
  A left multiplier~$L$ of $\CCINF(G)$ is of the form $f\mapsto D\ast
  f$ for a uniquely determined distribution $D\in\CCINF'(G)$.  A right
  multiplier is of the form $f\mapsto f\ast D$ for a uniquely
  determined distribution $D\in\CCINF'(G)$.  If a pair $(D_1,D_2)$ of
  distributions gives an element of $\Mult(A)$, then $D_1=D_2$.  Thus
  $\Mult(A)$ is the intersection of $\LMult(A)$ and $\RMult(A)$ inside
  $\CCINF'(G)$.
\end{lemma}

\begin{proof}
  Let $L\in\LMult(\CCINF(G))$.  Then we define a distribution
  $D_L\in\CCINF'(G)$ by $D_L(f)\defeq L(f)(1_G)$.  We view $\CCINF(G)$
  as an essential right module over $\CCINF(G)$ and~$L$ as a bounded
  module homomorphism.  The right module structure on $\CCINF(G)$ is
  the integrated form of the representation $\mu_G\cdot\rho$.
  Theorem~\ref{the:essential_smooth} yields that~$L$ is equivariant
  with respect to this representation of~$G$.  A straightforward
  computation now shows that $Lf = D_L\ast f$ for all $f\in\CCINF(G)$.
  If $D\ast f=0$ for all $f\in\CCINF(G)$, then $D\ast f(1)=0$ for
  all~$f$ and hence $D=0$.  Thus the distribution and the left
  multiplier $D\ast\blank$ determine each other uniquely.  The
  antipode on $\CCINF(G)$ extends to an algebra isomorphism between
  $\LMult(\CCINF(G))$ and $\RMult(\CCINF(G))$.  Hence the description
  of left multipliers above yields a description of right multipliers.
  If the pair $(D_1,D_2)$ determines a two-sided multiplier, then
  $(a\ast D_2)\ast b=a\ast (D_1\ast b)$ for all $a,b\in\CCINF(G)$.
  Thus the right multiplier associated to the distribution
  $(D_2-D_1)\ast b$ vanishes for all~$b$.  This implies $(D_2-D_1)\ast
  b=0$.  Since~$b$ is arbitrary, we obtain $D_2=D_1$.
\end{proof}

It remains to identify the distributions on~$G$ that give rise to
left, right and two-sided multipliers.  Let~$I$ be a
fundamental system of smooth compact subgroups of~$G$.  For $k\in I$
let~$\mu_k$ be the normalized Haar measure on~$k$, viewed as a
distribution on~$G$.  Thus the convolution with~$\mu_k$ on the left
and right averages a function over left or right $k$\nbd{}cosets.

\begin{proposition}  \label{pro:CCINF_mult}
  A distribution $D\in\CCINF'(G)$ is a left multiplier of $\CCINF(G)$
  if and only if $D\ast\mu_k\in\CINF'(G)$ for all $k\in I$ and a right
  multiplier if and only if $\mu_k\ast D\in\CINF'(G)$ for all $k\in
  I$.  There are bornological isomorphisms
  \begin{align*}
    \LMult(\CCINF(G)) &\cong
    \varprojlim_{k\in I} \CINF'(G/k)
    \cong
    \bigl(\varinjlim_{k\in I} \CINF(G/k)\bigr)';
    \\
    \RMult(\CCINF(G)) &\cong
    \varprojlim_{k\in I} \CINF'(k\backslash G)
    \cong
    \bigl(\varinjlim_{k\in I} \CINF(k\backslash G)\bigr)'.
  \end{align*}
\end{proposition}

\begin{proof}
  We only prove the isomorphisms for $\LMult(\CCINF(G))$.  The
  structure maps in the projective system $\CINF'(G/k)$ are right
  convolution with~$\mu_k$.  Recall that $\CCINF(G)=\varinjlim
  \CCINF(k\backslash G)$ and that left convolution with~$\mu_k$ is a
  projection onto $\CCINF(k\backslash G)$.  Thus
  $D\in\LMult(\CCINF(G))$ if and only if left convolution with
  $D\ast\mu_k$ is a bounded map from $\CCINF(k\backslash G)$ to
  $\CCINF(G)$.  Clearly, this is the case if $D\ast\mu_k$ has compact
  support.  Conversely, if $D\ast\mu_k$ does not have compact support,
  then there exist functions $(\phi_n)_{n\in\N}$ in
  $\CCINF(k\backslash G)$ whose support is contained in a fixed
  compact subset $L\subseteq G$ for which $D\ast\mu_k\ast\phi_n$ does
  not have a common compact support.  Multiplying the
  functions~$\phi_n$ by appropriate scalars we can achieve that
  $\{\phi_n\}$ is a bounded subset of $\CCINF(k\backslash G)$.  By
  construction, $D\ast\{\phi_n\}$ is not a bounded subset of
  $\CCINF(G)$, so that~$D$ is not a left multiplier.  Thus
  $D\in\LMult(\CCINF(G))$ if and only if $D\ast\mu_k$ has compact
  support for all $k\in I$.  An analogous computation for a set
  $S\subseteq\CCINF'(G)$ of distributions shows that~$S$ is bounded in
  $\LMult(\CCINF(G))$ if and only if $S\ast\mu_k$ is bounded in
  $\CINF'(G/k)$ for all $k\in I$.  This proves the first isomorphism.
  The second one follows from the universal property of direct limits.
\end{proof}

\begin{corollary}  \label{cor:CCINF_mult_prolie}
  If~$G$ is a projective limit of Lie groups, then
  $$
  \LMult(\CCINF(G))=\RMult(\CCINF(G))=\Mult(\CCINF(G)).
  $$
  If~$G$ is a Lie group then all three multiplier algebras are
  equal to $\CINF'(G)$.
\end{corollary}

The spaces $\CINF(G/k)$ for $k\in I$ are nuclear Fréchet spaces and
hence reflexive.  We can rewrite the inductive limit $\varinjlim_{k\in
I} \CINF(G/k)$ as a direct sum.  If~$G$ is metrizable, this is quite
easy: choose~$I$ to be a sequence and notice that $\CINF(G/k_n)$ is a
retract of $\CINF(G/k_{n+1})$ for any $n\in\N$.  If~$G$ is not
metrizable, the assertion is still correct, but the proof is more
complicated.  Therefore, $\varinjlim \CINF(G/k)$ is reflexive, so that
$\LMult(\CCINF(G))'\cong\varinjlim \CINF(G/k)$.  Furthermore, if~$G$
is countable at infinity, then Proposition~\ref{pro:CINF_smoothen}
shows that $\varinjlim_{k\in I} \CINF(G/k)$ is the smoothening of the
right regular representation on $\CINF(G)$.

\begin{proposition}  \label{pro:central_mult}
  Let $D\in\CCINF'(G)$.  Then~$D$ is a central multiplier of
  $\CCINF(G)$ if and only if $\mu_k\ast D\ast\mu_k\in
  \Center\CINF'(G//k)$ for all $k\in I$.  There is a natural
  isomorphism of bornological algebras
  $$
  \Center \Mult(\CCINF(G)) \cong \varprojlim \Center \CINF'(G//k).
  $$
\end{proposition}

\begin{proof}
  If~$D$ is a central multiplier of $\CCINF(G)$, then $\mu_k\ast
  D\ast\mu_k$ belongs to the center of
  $\mu_k\Mult(\CINF(G))\mu_k$.  Proposition~\ref{pro:CCINF_mult}
  yields an isomorphism of bornological algebras
  $\mu_k\Mult(\CINF(G))\mu_k=\CINF'(G//k)$.  Hence we have a bounded
  homomorphism $\Center\Mult(\CCINF(G))\to\varprojlim
  \Center\CINF'(G//k)$.

  Suppose conversely that $\mu_kD\mu_k$ be a
  central element of $\CINF'(G//k)$ for all $k\in I$.
  For any $j\in I$, $j\subseteq k$, $f\in\CCINF(G//k)$, we have
  $$
  \mu_j\ast D\ast f=
  \mu_j\ast D\ast\mu_j\ast f\ast\mu_k=
  f\ast\mu_j\ast D\ast\mu_j\ast \mu_k=
  f\ast\mu_k\ast D\ast\mu_k.
  $$
  Since this is independent of~$j$, we obtain $D\ast f=f\ast\mu_k\ast
  D\ast\mu_k$.  In particular, $D$ is a left multiplier.  A similar
  computation for $f\ast D$ shows $f\ast D=D\ast f$ because
  $\mu_kD\mu_k$ commutes with~$f$.  Hence~$D$ is central, so that we
  obtain an isomorphism $\Center \Mult(\CCINF(G)) \cong \varprojlim
  \Center \CINF'(G//k)$.  It is easy to check that it is bornological.
\end{proof}

If~$G$ is totally disconnected, then the spaces $G//k$ are all
discrete, so that $\CINF'(G//k)=\CCINF(G//k)$.  This special case is
covered in~\cite{Bernstein:Center}.  Now let~$G$ be a connected Lie
group.  If $[D,X]=0$ for all $X\in\LG$, then $[D,\delta_g]=0$ for all
$g\in G$ and hence~$D$ is central.  Thus a distribution is central if
and only if it commutes with~$\LG$.  In particular, the center of the
universal enveloping algebra of~$G$ is contained in the center of
$\CINF'(G)$.  The latter can be bigger than $Z\UniEnvel(G)$.  This
happens, for instance, if~$G$ has non-trivial center or if~$G$ is
compact.  However, there are also many Lie groups for which we have
$Z\UniEnvel(G)=Z\CINF'(G)$, that is, any central distribution is
supported at~$1$.  The following proposition only gives one class of
examples.

\begin{proposition}  \label{pro:center_complex_Lie}
  Let~$G$ be a connected complex Lie group with trivial center.  Then
  $\Center \Mult(\CCINF(G))$ is equal to the center of the
  universal enveloping algebra.
\end{proposition}

\begin{proof}
  Since~$G$ has trivial center, the adjoint representation of~$G$ on
  its Lie algebra~$\LG$ is faithful, so that $G\subseteq \Gl(\LG)$.
  Let $D\in\Center \CINF'(G)$ and $y\in\supp D$.  Since $\supp D$ is
  compact and conjugation invariant, the holomorphic function
  $$
  \C\ni s\mapsto \exp(sX)y\exp(-sX) \in \Gl(\LG)
  $$
  is bounded for any $X\in\LG$.  Liouville's Theorem yields that it is
  constant, that is, $[X,y]=0$.  This implies $\supp D=\{1\}$
  because~$G$ has trivial center.  Now use the identification of
  distributions supported at~$1$ with the universal enveloping
  algebra.  Since~$G$ is connected, a distribution is central if and
  only if it commutes with~$\LG$.
\end{proof}

\end{document}